\documentclass[12pt]{birkjour}
\usepackage{amssymb,amsfonts,amsmath}

\renewcommand{\theequation}{\thesection.\arabic{equation}}
\newtheorem{Theorem}{Theorem}[section]
\newtheorem{Lemma}[Theorem]{Lemma}
\newtheorem{Definition}[Theorem]{Definition}
\newtheorem{Remark}[Theorem]{Remark}
\begin{document}

\title[Pullback dynamics of a 3D Navier-Stokes equation with nonlinear viscosity] {Pullback dynamics of a 3D Navier-Stokes equation with nonlinear viscosity}
\author{Xin-Guang Yang}
\address{Department of Mathematics and Information Science, Henan Normal University, Xinxiang 453007, China}
\email{yangxinguang@hotmail.com}

\author{Baowei Feng}
\address{College of Economic Mathematics, Southwestern University of Finance and Economics, Chengdu 611130, China}
\email{bwfeng@swufe.edu.cn}

\author{Shubin Wang}
 \address{School of Mathematics and Statistics, Zhengzhou University, Zhengzhou 450001, China}
 \email{wangshubin@zzu.edu.cn}

\author{Yongjin Lu} %\corref{cor1}}
\address{Department of Mathematics and Economics, Virginia State University, Petersburg, VA 23806, USA.}
\email{ylu@vsu.edu}

\author{To Fu Ma (corresponding)}
 \address{Institute of Mathematical and Computer Sciences, University of S\~ao Paulo, 13560-970 S\~ao Carlos, SP, Brazil}
 \email{matofu@icmc.usp.br}


\begin{abstract}
		This paper is concerned with pullback dynamics of 3D Navier-Stokes equations with variable viscosity and subject to time-dependent external forces.
		Our main result establishes the existence of finite-dimensional pullback attractors in a general setting involving tempered universes.
		We also present a sufficient condition on the viscosity coefficients that guarantees the attractors are nontrivial.
		We end the paper by showing the upper semi-continuity of pullback attractors as the non-autonomous perturbation vanishes.
\end{abstract}

\maketitle

%\begin{keyword}
%Navier-Stokes equations, nonlinear viscosity, pullback attractors, fractal dimension, Ladyzhenskaya model.

%\MSC[2010] 35B40 \sep 35B41 \sep 35Q30 \sep 76D03 \sep 76D05.
%\end{keyword}

% \linenumbers

\section{Introduction}
\setcounter{equation}{0}

The incompressible 3D Navier-Stokes equations
\begin{equation}\label{bmy1-a}
\begin{cases}
u_t-\nu\Delta u + (u\cdot \nabla)u+\nabla p = f,  & \\
\nabla \cdot \, u=0, &
\end{cases}
\end{equation}
defined in a domain $\Omega$ of $\mathbb{R}^3$, is a classical mathematical model for viscous fluids with several applications in real world problems. The unknown variables $u,p$ stand for, respectively, the velocity field and the scalar pressure of the fluid, and parameter $\nu >0$ is a kinematic viscous coefficient. Its mathematical analysis has been studied by many authors since the pioneering works of Leray and Hopf \cite{leray1933, leray1934, leray193401,hopf1951}.

Because the study of Navier-Stokes equations is very important in the understanding of fluid turbulence, Ladyzhenskaya
(in the 1960s) asked if the system \eqref{bmy1-a} could determine completely the motion of a viscous fluid flow.
In fact, the rigorous answer to the well-posedness of \eqref{bmy1-a}, with given boundary condition
$u|_{\partial\Omega}=a$ and appropriate initial data $u|_{t=0}=b(x)$, is still an open problem
(see e.g. \cite{lions1998}). However, in certain cases, adding a physically meaningful feature may provide global existence.
For instance, some variants of \eqref{bmy1-a} obtained by adding strong viscosity or time-delay,
are well-posed. Then Ladyzhenskaya \cite{lady,L2} proposed a class of
Navier-Stokes equations of the form
\begin{equation}\label{bmy1-b}
\begin{cases}
u_t- \mbox{div}\Big[ \big(\nu_0+\nu_1 \|\mathcal{D}u\|^2_{L^2} \big) \mathcal{D}u\Big] + (u\cdot \nabla)u+\nabla p = f(x,t), & \\
\nabla \cdot \, u=0,&
\end{cases}
\end{equation}
where $\mathcal{D}u=\nabla u+\nabla u^{\top}$.

Due to the variable viscosity, the system (\ref{bmy1-b}) has many desired features, namely,
meaningful physical justification, well-posedness under appropriate initial-boundary conditions in some sub-critical cases,
and satisfies the Stokes principle. We refer the reader to Beir\~ao da Veiga \cite{bei} for further modeling aspects and a
regularity result.

On the other, the uniqueness and stability of solutions when Reynold number is large are open questions.
To overcome these difficulties and simplify a little (\ref{bmy1-b}), Lions \cite{JLL1} considered
the model with nonlinear viscosity,
\begin{equation}\label{bmy1-1}
\begin{cases}
u_t - (\nu+\nu_0\| \nabla u \|_{L^2}^2)\Delta u + (u\cdot \nabla)u + \nabla p = f(x,t), \;\ x\in\Omega, \ t > \tau,& \\
\nabla \cdot \, u=0,&\\
u|_{\partial\Omega}= 0,&\\
u(x,\tau) = u_{0}(x),\ x\in\Omega, &
\end{cases}
\end{equation}
where $\Omega$ is a bounded domain of $\mathbb{R}^3$ with smooth boundary $\partial \Omega$
and $\nu,\nu_0>0$. Note that for $\nu_0 = 0$ this system reduces formally to \eqref{bmy1-a}.

The continuity between the classical Navier-Stokes equations \eqref{bmy1-a} and the approximating system \eqref{bmy1-1} is still an open problem. In particular, the convergence of the pullback attractors of \eqref{bmy1-1} to trajectory attractors of \eqref{bmy1-a} is also unknown and seems to be an interesting problem.

Next, we recall some particularly important results pertaining to \eqref{bmy1-1}. The existence of weak global solutions and uniqueness
were proved in Lions \cite[Chapter 2]{JLL1}. With respect to non-autonomous dynamics, the existence of uniform attractors was studied in \cite{cys}.
However, as far as we know, there are no results on pullback dynamics to the system \eqref{bmy1-1} with time-dependent external forces.

The main results of the present paper are summarized as follows:

$(i)$ We establish the existence of minimal (unique) pullback attractors of \eqref{bmy1-1} and their structure in a general framework featuring the notion of tempered universe. The relation between families of pullback attractors in various universes is also obtained. See Theorems \ref{th5.1} and \ref{re}.

$(ii)$ We provide estimates on the (finite) fractal dimension of pullback attractors in $H$.
Our result differs from those for 2D and 3D classical Navier-Stokes equations since here the variable viscosity coefficient
plays a special role. See Theorem \ref{th5-f}.

$(iii)$ We present a sufficient condition ensuring that our pullback attractors are nontrivial.
See Theorem \ref{th5-g}.

$(iv)$ The upper semi-continuity of pullback attractors as the norm of external force vanishes
is also studied. See Theorem \ref{th5.10}.

\medskip

The rest of this paper is arranged as follows. We present some preliminaries and notations in Section \ref{sec2}.
In Section \ref{sec-3} we present our main results and make some remarks.
In Section \ref{sec5} we complete the proofs of main results.

\section{Preliminaries}\label{sec2}
\setcounter{equation}{0}

\subsection{Notations and functional spaces}
We shall work with usual Sobolev spaces $H^s_0(\Omega)$ and $H^s(\Omega) \; (s\in \mathbb{R})$.
Let us define
$$
E =\{u \in (C^{\infty}_0(\Omega))^3 \, | \; \nabla\cdot u=0\} \quad \mbox{and} \quad H = \overline{E},
$$
with closure in $(L^2(\Omega))^3$ topology. The inner product and
norm of $H$ are denoted by $(\cdot,\cdot)$ and $|\cdot|$, respectively, and defined by
\begin{align*}
(u,v) & =\sum^3_{j=1}\int_{\Omega}u_j(x)v_j(x)d x, \quad \forall\ u, v\in 	 (L^2(\Omega))^3,\\
|u|^2 &= (u,u),  \quad \forall\ u\in (L^2(\Omega))^3.
\end{align*}
We denote by $V$ the closure of $E$ in $(H^1(\Omega))^3$ topology. The inner product and norm in $V$
are denoted by $((\cdot,\cdot))$ and $\|\cdot\|$, respectively, and defined by
\begin{align*}
((u,v)) & =\sum^3_{i,j=1}\int_{\Omega}\frac{\partial u_j}{\partial x_i}\frac{\partial v_j}{\partial x_i}d x,\quad \forall\ u, v\in (H^1_0(\Omega))^3,\\
\|u\|^2 & =((u,u)), \quad \forall\ u\in (H^1_0(\Omega))^3.
\end{align*}
The dual spaces of $H$ and $V$ are denoted by $H'$ and $V'$, respectively, and the injections $V\hookrightarrow H\equiv H'\hookrightarrow V'$ are
dense and continuous. The usual notations $\|\cdot\|_{*}$ and $\langle \cdot , \cdot \rangle$ stand for the
norm in $V'$ and the duality pairing between $V$ and $V'$ (also $H$ and its dual space), respectively.

Let $X$ be a Banach space with norm $\|\cdot\|_{X}$. The Hausdorff semi-distance $\mbox{dist}_{X}$ between two sets $B_{1},B_{2}\subseteq X$ is defined by
$$
{\rm dist}_{X}(B_{1},B_{2})=\sup\limits_{x\in B_{1}}\inf\limits_{y\in B_{2}} \parallel x-y \parallel_{X}.
$$
This will be used in the definition of attractors.

\subsection{Abstract equivalent equation}
Let $P$ be the Helmholz-Leray orthogonal projection operator  from $(L^2(\Omega))^3$ onto $H$.
We define $A =-P\Delta $, the Stokes operator,
with domain $D(A)=(H^{2}(\Omega))^3\bigcap(H^{1}_{0}(\Omega))^3$. Then the operator $A: V\rightarrow V'$ has the property $\langle Au,v\rangle=((u,v))$,
for all $u, v\in V$, and it is an isomorphism from $V$ into $V'$.
The eigenvalues of $A$ with Dirichlet boundary condition in $(L^2(\Omega))^3$ are denoted by $\{\lambda_j\}^{\infty}_{j=1}$ (and satisfy
$0<\lambda_1\leq\lambda_2\leq\cdots$). The corresponding eigenfunctions are denoted by $\{\omega_j\}^{\infty}_{j=1}$ and provide an orthogonal
basis for $H$.

We define the fractional operator $A^{s}$ $(s\in \mathbb{C})$ (see \cite{boyer,te}) as follows.
\begin{align}
A^{s}u & =\sum_{i=1}^{\infty}\lambda_{j}^{s}(u,\omega_j)\omega_{j}, \; u \in H, \; s \in \mathbb{C}, \; j\in \mathbb{N}, \label{2}\\
V^s & = D(A^{s})=\Big\{u\in H: A^{s}u\in H,\ \sum_{i=1}^{\infty}\lambda_i^{2 s}|(u,\omega_i)|^2<+\infty\Big\}, \nonumber \\
|A^{s}u | & =\Big(\sum^{\infty}_{i=1}\lambda_i^{2 s}|(u,\omega_i)|^2\Big)^{1/2}, \nonumber
\end{align}
where $D(A^{s})$ denotes the domain of $A^{s}$ equipped with the inner product and the norm
$$
(u,v)_{V^s}=(A^{\frac{s}{2}}u, A^{\frac{s}{2}}v),\quad\ \|u\|_{V^s}^2=(u,u)_{V^s}.
$$
In particular, $V=V^1$ and $V^{2}= W=(H^{2}(\Omega))^3\bigcap (H^{1}_{0}(\Omega))^3$.
%For convenience, we often use the norm notation $\|\cdot\|_{V^{\frac{s+1}{2}}}=\|\cdot\|_{1+s}$.
In addition, the immersion
$$
D(A^{\frac{s}{2}})\hookrightarrow D(A^{\frac{r}{2}}), \quad s > r,
$$
is continuous and
$$
D(A^{\frac{s}{2}})\hookrightarrow\hookrightarrow
(L^{\frac{6}{(3-2s)}}(\Omega))^{3}, \quad \mbox{$0 < s \le \frac{3}{2}$},
$$
is compact. The abstract setting of \eqref{bmy1-1} also uses an operator $\mathbb{A}:V\to V'$ defined by $\mathbb{A}u=-\nu_0\|u\|^2\Delta u$,
that satisfies
$$
\langle\mathbb{A}u,v\rangle =\nu_0\|u\|^2\langle -\Delta u,v\rangle= \nu_0\|u\|^2((u,v)),\ \forall\ u,v\in V.
$$
Noting that $\mathbb{A}$ is monotone from $V$ into $V'$, we have
\begin{eqnarray}
\|\mathbb{A}u\|_{\mathcal{L}(V,V')}=\sup_{\|v\|=1,\ v\in V}|\langle \mathbb{A}u,v\rangle|=\sup_{\|v\|=1,\ v\in V} \nu_0\|u\|^2|a(u,v)|\leqslant \nu_0\|u\|^3.\nonumber
\end{eqnarray}
Finally we define the bilinear and trilinear operators (see \cite{te}),
\begin{align*}
B(u,v) & = P((u\cdot\nabla)v),  \\
b(u,v,w) & =(B(u,v),w)=\sum^3_{i,j=1}\int_{\Omega}u_i\frac{\partial v_j}{\partial x_i} w_jd x ,
\end{align*}
for all $u,v,w \in E$.

\section{Main results and comments}\label{sec-3}
\setcounter{equation}{0}

\subsection{Well-posedness and regularity of global solutions}

A function $u=u(x,t)$ in $L^{\infty}(\tau,T;H)\cap L^{4}(\tau,T;V)$ is called weak solution of \eqref{bmy1-1} if
$$
\begin{cases}
\frac{d}{dt} (u(t),v)+(\nu+\nu_0\|\nabla u(t)\|^2)((u(t),v))+b(u(t),u(t),v)=\langle f(x,t),v\rangle,   \\
u(\tau)=u_{0},
\end{cases}
$$
for all $v\in V$, in the sense of distributions. Based on previous notations, the non-autonomous system \eqref{bmy1-1}
can be rewritten as an abstract functional equation
\begin{equation}\label{yang-2}
\begin{cases}  u_t + \nu Au+\mathbb{A}u+B(u,u)=f(x,t), \; t \ge \tau,  & \\
\nabla\cdot u=0, & \\
u|_{\partial\Omega}=0,& \\
u(x,\tau)=u_{0}(x).&
\end{cases} \end{equation}

\begin{Theorem} [Well-posedness] \label{thm1} Assume the external force $f=f(x,t) \in L_{loc}^{4/3}(\mathbb{R},V')$ and $u_0\in H$. Then equation \eqref{yang-2} has a unique weak solution
	$$
	u \in C(\tau,T;H)\cap L^4(\tau,T; V),
	$$
	that depends continuously on initial data.
\end{Theorem}

\paragraph{Proof} The proof of global existence and uniqueness of weak solutions is similar to the ones in
\cite{JLL1,JLLandP} by using Faedo-Galerkin approximations. Here we omit the details.
Since the solution is continuously dependent on the initial data and $\frac{du}{dt}\in L^2(\tau,T;V')$, from  Aubin-Lions lemma we derive that $u \in C(\tau,T;H)$. \qed

\begin{Remark} We recall that an evolution process on a metric space $X$ is a two-parameter family of mappings  $S(t,\tau)$, $t \ge \tau$,
	defined on $X$ such that: $(i)$ $S(t,t) = Id_X$ for all $t\in \mathbb{R}$ (identity).
	$(ii)$ $S(t,\tau)=S(t,s)S(s,\tau)$, $t \ge s \ge \tau$ (semigroup property).
	Moreover, a process is called continuous if $(t,\tau,z) \mapsto S(t,\tau)z$ is continuous for $t \ge \tau$, $z\in X$.
	Then Theorem \ref{thm1} shows that weak solutions of problem \eqref{yang-2} generates a continuous evolution process in $H$.
\end{Remark}

Regarding to regular solutions following result holds.

\begin{Theorem}[Regularity] \label{thm2} If $u_0\in D(A^{\frac{\sigma}{2}})$ with
	$\sigma \in [0,1]$ and $f \in L^2_{loc}(\mathbb{R};H)$, then the solution of \eqref{yang-2} satisfies
	$$
	u\in C(\tau,T;D(A^{\frac{\sigma}{2}}))\cap L^2(\tau,T;D(A^{\frac{\sigma+1}{2}})).
	$$
	Moreover, the regular solutions generates a continuous evolution process $S(t,\tau)$ on $D(A^{\frac{\sigma}{2}})$.
\end{Theorem}

\paragraph{Proof} The arguments are similar to the ones in \cite{JLL1,JLLandP} and the proof of Lemma \ref{le4.8}. \qed

\subsection{Forward dynamical systems} \label{for}

With respect to autonomous dynamics, when $f=f(x)$ is time independent, we have:

\begin{Theorem} [Global attractors] \label{re4.6}
	Assume that $f \in V'$. Then, with respect to weak solutions, problem \eqref{bmy1-1} has a global
	attractor $\mathcal{A}^H \in H$. Moreover, with respect to regular solutions, for $f \in H$, problem \eqref{bmy1-1} has a global attractor $\mathcal{A}^{D(A^{\frac{1}{4}})} \in D(A^{\frac{1}{4}})$.
\end{Theorem}

The proof of Theorem \ref{re4.6} is given in Section \ref{sec4}.

\begin{Remark}
	In the non-autonomous case where $f$ is an external force in the class of
	{\it translation bounded}, {\it translation compact} or {\it normal} functions (see \cite{cv3}),
	the existence of uniform attractors was established in \cite{cys}. \qed
\end{Remark}

\subsection{Pullback dynamics}\label{pul}

In order to present our results on pullback dynamics of problem \eqref{bmy1-1} we recall some basic definitions of the theory (cf. \cite{clr2006,clr,gmr,luk}).

\begin{Definition} \label{def-universe} A universe $\mathcal{D}$ defined in a metric space $X$ is a class of families $\hat{D}$ of the form
	$\hat{D}=\{D(t)|t\in\mathbb{R}\}$, where each $D(t)$ is a nonempty bounded subset of $X$.
\end{Definition}

\begin{Definition} Given a universe $\mathcal{D}$ defined on $X$,
	a family $\mathcal{A} = \{\mathcal{A}(t)\}_{t \in \mathbb{R}}$ is a called pullback $\mathcal{D}$-attractor
	of a process $S(t,\tau):X \to X$ if the following properties hold:
	\begin{itemize}
		\item[$(i)$] Compactness: $\mathcal{A}(t)$ is a nonempty compact set of $X$, $t \in \mathbb{R}$,
		\item[$(ii)$] Invariance: $U(t,\tau)\mathcal{A}(\tau)=\mathcal{A}(t)$, $t \ge \tau$ ,
		\item[$(iii)$] Pullback $\mathcal{D}$-attraction: $\displaystyle \lim_{\tau \to -\infty}{\rm dist}_{X}(S(t,\tau)D(\tau),\mathcal{A}(t))=0$,
		$t\in \mathbb{R}$, $\{D(t)\}_{t\in\mathbb{R}} \in \mathcal{D}$.
	\end{itemize}
	In addition, a $\mathcal{D}$-attractor $\mathcal{A}$ is said to be minimal if whenever $\hat{C}$ is another
	$\mathcal{D}$-attracting family of closed sets, then $\mathcal{A}(t) \subset C(t)$, for all $t \in \mathbb{R}$.
\end{Definition}

Our results are concerned with families of universes determined by the time-dependent force $f=f_{\mu}$.
Let us put $\mu_0=\nu\lambda_1$ and assume either
\begin{eqnarray}
f \in L^2_{loc}(\mathbb{R};V') \; \mbox{and} \; \int^{t}_{-\infty}e^{\mu s}\|f(s)\|^2_{V'} \,ds
< +\infty \;\; \mbox{for some} \;\; \mu \in (0,\mu_0], \; t \in\mathbb{R},  \label{f7a}
\end{eqnarray}
or
\begin{eqnarray}
f \in L^2_{loc}(\mathbb{R};H) \; \mbox{and} \; \int^{t}_{-\infty}e^{\mu s} |f(s)|^2 \, ds
< +\infty \; \mbox{for some} \;\; \mu \in (0,\mu_0], \; t\in\mathbb{R}.  \label{f7a-1}
\end{eqnarray}
The corresponding $\mu$-indexed universes are defined by
\begin{equation} \label{gg7}
\mathcal{D}^H_{\mu} = \left\{ \hat{D}_{\mu} \, | \, D_{\mu}(t) \subset B(0,\rho_{\hat{D}_{\mu}}(t)) \; \mbox{with} \;
\lim_{\tau \to -\infty}	|\rho_{\hat{D}_{\mu}}(\tau)|^2e^{\mu\tau} = 0 \right\},
\end{equation}
where $\rho_{\hat{D}_{\mu}}: \mathbb{R} \to \mathbb{R}^{+}$ is a continuous function.
Then, $B(0,\rho_{\hat{D}_{\mu}}(t))$ denotes a family of closed ball in $H$ of radius $\rho_{\hat{D}_{\mu}}(t)$.
Analogously,
we define
\begin{equation} \label{gg7-1}
\mathcal{D}^{D(A^{\frac{1}{4}})}_{\mu} = \left\{ \hat{D}_{\mu} \, | \, D_{\mu}(t) \subset \tilde{B}(0,\rho_{\hat{D}_{\mu}}(t)) \; \mbox{with} \;
\lim_{\tau \to -\infty}	|\rho_{\hat{D}_{\mu}}(\tau)|^2e^{\mu\tau} = 0 \right\},
\end{equation}
where $\tilde{B}(0,\rho_{\hat{D}_{\mu}}(t))$ are closed balls in
$D(A^{\frac{1}{4}})$. We note that these universes are inclusion closed.

In this direction our first result is:

\begin{Theorem} [Pullback attractors] \label{th5.1}
	Under assumption \eqref{f7a}, the process generated by the weak solutions of problem \eqref{bmy1-1} possesses
	a minimal family of pullback $\mathcal{D}_{\mu}^{H}$-attractors $\mathcal{A}_{\mu}^{H}$ in $H$, for all $\mu \in (0,\mu_0]$. Analogously, under assumption \eqref{f7a-1}, the process generated by the regular solutions of problem \eqref{bmy1-1} possesses a minimal family of pullback $\mathcal{D}_{\mu}^{D(A^{\frac{1}{4}})}$-attractors $\mathcal{A}^{D(A^{\frac{1}{4}})}_{\mu}$ in $D(A^{\frac{1}{4}})$, for all $\mu \in (0,\mu_0]$.
\end{Theorem}

The proof of Theorem \ref{th5.1} is given in Section \ref{se3.2}.

\begin{Remark} Given $\hat{D} = \{D(t)\}_{t \in \mathbb{R}}$ in $\mathcal{D}_{\mu}^H$, the pullback omega limit of $\hat{D}$, at time $t$,
	is defined by
	$$
	\Lambda(\hat{D},t) =
	\bigcap_{\tau \leq t} \, \overline{\bigcup_{s\leq \tau} S(t,s) D(s)}^H.
	$$
	Then, from an abstract existence result (e.g. \cite[Theorem 3.11]{gmr}), the pullback attractor $\mathcal{A}_{\mu}^{H}=\{\mathcal{A}_{\mu}^{H}(t)\}_{t \in \mathbb{R}}$ in Theorem \ref{th5.1} is characterized by	
	$$
	\mathcal{A}_{\mu}^{H}(t) = \overline{\bigcup_{\hat{D} \in \mathcal{D}_{\mu}^{H}} \Lambda(\hat{D},t)}^H, \quad  t \in \mathbb{R}.
	$$
	Analogous characterization holds for $\mathcal{A}^{D(A^{\frac{1}{4}})}_{\mu}$. \qed  	
\end{Remark}

\subsection{Fractal dimension of attractors}

Let $K$ be a non-empty compact subset of a Hilbert space $H$. Given $\varepsilon>0$, we denote $N_{\varepsilon}(K)$ to be the minimum number of open balls in $X$ with radius $\varepsilon$ which are necessary to cover $K$. Then the fractal dimension of $K$ is defined by
$$
dim_F(K)=\displaystyle{\limsup_{\varepsilon\rightarrow 0^+}} \, \frac{\log(N_{\varepsilon}(K))}{\log(\frac{1}{\varepsilon})}.
$$
We say that a pullback attractor $\mathcal{A}=\{A(t)\}_{t\in\mathbb{R}}$ has finite fractal dimension if there exists $d>0$ such that
$dim_F(A(t)) \le d$ for all $t\in \mathbb{R}$. We find upper bounds of fractal dimension for pullback attractors in $H$ by verifying the uniform differentiability of the process and using a trace formula \cite{clr}.

For the next result we define
\begin{equation} \label{G}
\langle h \rangle|_{\leq t} = \limsup_{s \rightarrow -\infty} \frac{1}{t-s} \int^t_s h(r) dr
\quad \mbox{and} \quad
{\rm G}(t) = \frac{\| f \|^2_{L^{\infty}(-\infty,t; H)}}{\nu^2_0\lambda_1},
\end{equation}
where ${\rm G}(t)$ is the Grashof number for the non-autonomous system \eqref{bmy1-1}. In the autonomous case, one has
${\rm G}_{a} = \frac{|f|^2}{\nu^2_0\lambda_1}$.

\begin{Theorem} [fractal dimension] \label{th5-f} The pullback attractors $\mathcal{A}_{\mu}^{H}$ given by Theorem \ref{th5.1}
possess finite fractal dimension, for any $\mu \in (0 , \mu_0]$. Moreover, we have the estimate
\begin{equation} \label{FD}
	dim_F(\mathcal{A}^H_{\mu}(t))\leq \max\Big\{3, C_F {\rm G}(t) + \frac{2\nu}{27}\Big\},
\end{equation}
where $C_F>0$ is a constant.
\end{Theorem}	

The proof of Theorem \ref{th5-f} is given in Section \ref{sec4.3}. There in Remark \ref{rem-fractal}, we provide  further estimates for \eqref{FD} and some comments on how variable viscosity in \eqref{bmy1-1} plays an important.

\bigskip

Next we present a criterion that guarantees the pullback attractors are nontrivial, that is, not a one point set.
We shall use a generalized Grashof number
\begin{equation} \label{Gg}
{\rm G^g}(t) = \frac{ \langle |f|^2 \rangle|_{\leq t}^{\frac{1}{2}} }{\nu^2_0\lambda_1}.
\end{equation}

\begin{Theorem} \label{th5-g} If
$$
{\rm G^g}(t) \geq \sqrt{\frac{\nu_0}{c\nu+4\nu_0^2\nu\lambda_1}},
$$
	then the pullback attractor $\mathcal{A}^H_{\mu}$ is not a single-trajectory set.
\end{Theorem}

The proof of Theorem \ref{th5-g} is given in Section \ref{se3.8}. \qed

\subsection{Comparing families of pullback attractors}

In the literature, several studies were concerned with pullback attractors for evolution systems with uniformly bounded absorbing family, e.g. \cite{clr2006,scd06,wzz}. To problem \eqref{bmy1-1},
as we show in Lemma \ref{le5.3}, our pullback $\mathcal{D}_{\mu}$-absorbing family needs not to be uniformly bounded, which is a key difference.

Apart the universes defined in \eqref{gg7}-\eqref{gg7-1},
we also consider $\mathcal{D}_F^H$ and $\mathcal{D}_F^{D(A^{\frac{1}{4}})}$, the universe of fixed nonempty bounded subsets of $H$ and $D(A^{\frac{1}{4}})$, respectively. That is,
$$
\mathcal{D}_F^X = \{\hat{D} \, | \, D(t)=D, \; \forall \, t \in \mathbb{R}, \mbox{where $D$ is a bounded set of $X$} \}.
$$
Then we see that $\mathcal{D}_F^H$ and $\mathcal{D}_F^{D(A^{\frac{1}{4}})}$ satisfy
\eqref{gg7}-\eqref{gg7-1}, and
$$
\mathcal{D}_F^H \subset \mathcal{D}^H_{\mu} \subset \mathcal{D}^H_{\mu_0} \quad \mbox{and} \quad
\mathcal{D}_F^{D(A^{\frac{1}{4}})} \subset \mathcal{D}_{\mu}^{D(A^{\frac{1}{4}})} \subset \mathcal{D}_{\mu_0}^{D(A^{\frac{1}{4}})}.
$$
However, $\mathcal{D}_F^H$ and $\mathcal{D}_F^{D(A^{\frac{1}{4}})}$ are not inclusion closed.
To establish our result we denote by $\mathcal{A}_{F}^H$ and $\mathcal{A}_{F}^{D(A^{\frac{1}{4}})}$ the corresponding pullback attractors given by Theorem \ref{th5.1}. We also take
$\hat{D}_{\mu_0} = \{D_{\mu_0}(t)\}_{t\in \mathbb{R}}$ the $\mathcal{D}_{\mu_0}$-absorbing family defined in $H$.

\begin{Theorem}\label{re} In the context of Theorem \ref{th5.1} we have:
\begin{enumerate}
	\item $\mathcal{A}^H_F(t)\subset \mathcal{A}^H_{\mu}(t)\subset \mathcal{A}^H_{\mu_0}(t)$, \quad $0< \mu \leq \mu_0$, \quad $t \in \mathbb{R}$.
	
	\item $\mathcal{A}^{D(A^{\frac{1}{4}})}_F(t)\subset \mathcal{A}^{D(A^{\frac{1}{4}})}_{\mu}(t)\subset \mathcal{A}^{D(A^{\frac{1}{4}})}_{\mu_0}(t)$, \quad $0< \mu \leq \mu_0$, \quad $t \in \mathbb{R}$.
	
	\item If $\displaystyle \bigcup_{t \le T} D_{\mu_0}(t)$ is bounded, then $\mathcal{A}^H_F(t)= \mathcal{A}^H_{\mu}(t)= \mathcal{A}^H_{\mu_0}(t)$, \quad $t\le T$.
	
	\item If $\displaystyle \bigcup_{t \le T} \mathcal{D}^{D(A^{\frac{1}{4}})}_{\mu_0}$ is bounded, then
	$\mathcal{A}^{D(A^{\frac{1}{4}})}_F(t)= \mathcal{A}^{D(A^{\frac{1}{4}})}_{\mu}(t)
	= \mathcal{A}^{D(A^{\frac{1}{4}})}_{\mu_0}(t)$, \quad $t\le T$.
\end{enumerate}	
\end{Theorem}

The proof of Theorem \ref{re} is presented in Section \ref{pro}.

\subsection{Upper semi-continuity of attractors}
Let us consider the perturbed system \eqref{bmy1-1} with $f(x,t) = \varepsilon h(x,t)$, for some $\varepsilon>0$ and $h \in L^2_{loc}(\mathbb{R},H)$ satisfying \eqref{f7a-1}. Then from Theorem \ref{th5.1}, with respect to the universe $\mathcal{D}_F$,
the dynamics of system \eqref{bmy1-1} possesses a pullback attractor
$\mathcal{A}_{\varepsilon} = \{ \mathcal{A}_{\varepsilon}(t)\}_{t \in  \mathbb{R}}$ that attracts any bounded set of $H$.
Let us denote by $\mathcal{A}_0 \subset H$ the corresponding global attractor of \eqref{bmy1-1} with $f=0$.
The upper semi-continuity of $\mathcal{A}_{\varepsilon}(t)$ to ${\mathcal{A}}_0$ as $\varepsilon \to 0$ in $H$ is stated in the following theorem.

\begin{Theorem} [Upper semi-continuity] \label{th5.10}
In the above context, suppose that
\begin{equation} \label{ya-1}
\int_{-\infty}^t e^{\eta s}|h(s)|^2 \, ds < \infty, \quad \eta \in (0, \nu \lambda_1].
\end{equation}
Then the pullback attractor ${\mathcal{A}}_{\varepsilon}$ is upper semi-continuous
with respect to $t\to 0$, that is,
$$
\lim_{\varepsilon \to 0^{+}} {\rm dist}_{H}({\mathcal {A}}_{\varepsilon}(t), {\mathcal {A}_0}) = 0.
$$
for all $t\in \mathbb{R}$.
\end{Theorem}

The proof of Theorem \ref{th5.10} is presented in Section \ref{sub5.2} and it is based on arguments from \cite{wq2010}.

\section{Proof of main results}\label{sec5}
\setcounter{equation}{0}

In this section, we shall present the proof of main results.

\subsection{Theory for pullback dynamics}

Some basic idea on evolutions processes, universes and pullback attraction were present previously in Section \ref{sec-3}.

\begin{Definition}
	We say that a family ${\mathcal B}=\{B(t)\}_{t\in{\mathbb R}}$ is pullback absorbing for a process $\{U(t,\tau)\}$ on $X$, if for every $t\in \mathbb{R}$ and any bounded subset $B\subset X$, there exists a time $T(t,B) \le t$, such that $U(t,\tau)B\subset B(t)$ if $\tau \le T(t,B)$. In addition, given a universe $\mathcal{D}$, we say that ${\mathcal B}$ is a pullback $\mathcal{D}$-absorbing family if, for every $t \in \mathbb{R}$ and $\hat{D} \in \mathcal{D}$,
	there exists $T(t,\hat{D}) < t$ such that $U(t,\tau)D(\tau) \subset B(t)$ for $\tau \le T(t,\hat{D})$.
\end{Definition}

\begin{Definition}
	Let ${\mathcal B}=\{B(t)\}_{t\in{\mathbb R}}$ be a family of subsets of $X$. We say that a process $\{U(t,\tau)\}$ defined on $X$ is pullback $\mathcal{B}$-asymptotically compact, if for any sequences $\tau_{n}\rightarrow -\infty$ and $x_{n}\in B(\tau_{n})$,
	the sequence $\{U(t,\tau_{n})x_{n}\}$ is precompact in $X$ for all $t\in \mathbb{R}$.
\end{Definition}

In what follows we recall a useful compactness criterion, known as condition (C), proposed by Ma, Wang and Zhong \cite{mwz,wz}. Here we call it the (MWZ)-condition.

\begin{Definition} Let $X$ be a Banach space and $\mathcal{D}$ a given universe. We say that a process
	$\{U(t,\tau)\}$ defined on $X$ satisfies the (MWZ) condition (with respect to $\mathcal{D}$) if for any $\hat{B} = \{B(t) \}_{t\in\mathbb{R}} \in \mathcal{D}$, for any fixed  $t\in\mathbb{R}$, any $\varepsilon>0$, there exists a pullback time $\tau_{\varepsilon}=\tau(t,\varepsilon,\hat{B})\leq t$ and a finite dimensional subspace $X_1\subset X$ such that
\begin{itemize}
\item[$(i)$] $P\left( \bigcup_{s\leq \tau_{\varepsilon}}U(t,s)B(s)\right)$ is bounded,
\item[$(ii)$] $\| (I-P) \left( \bigcup_{s\leq \tau_{\varepsilon}}U(t,s)B(s) \right) \|_X \leq \varepsilon$,
\end{itemize}
where $P$ is the bounded projection from $X$ to $X_1$.
\end{Definition}

The proof of existence of pullback attractors for problem \eqref{bmy1-1} is based on the following result.

\begin{Theorem} {\bf (\cite[Theorem 3.11]{gmr})} \label{thm-EPA} Let $\{U(t,\tau)\}$ be a continuous evolution process defined on a Banach space $X$ and let $\mathcal{D}$ be a universe defined on $X$. Suppose that
\begin{itemize}
	\item[$(i)$] $\{U(t,\tau)\}$ satisfies (MWZ) with respect to $\mathcal{D}$;
	\item[$(ii)$] $\{U(t,\tau)\}$ admits a pullback $\mathcal{D}$-absorbing family $\hat{D}_0=\{D_0(t)\}_{t\in\mathbb{R}}$.
\end{itemize}
Then, $U(t,\tau)$ possesses a minimal pullback $\mathcal{D}$-attractor
$\mathcal{A}_{\mathcal{D}} = \{A_{\mathcal{D}}(t)\}_{t \in \mathbb{R}}$.
Moreover, if $\hat{D}_0\in\mathcal{D}$ and $\mathcal{D}$ is inclusion closed, then $\mathcal{A} \in \mathcal{D}$.
\end{Theorem}

\subsection{Proof of Theorem \ref{re4.6}: Global attractor for autonomous case $f=f(x)$} \label{sec4}
In this subsection, we will prove the existence, regularity of global attractor for autonomous case of \eqref{bmy1-1}.

\begin{Lemma} [\cite{leray1934,te}]
	The bilinear operator $B(u,v)$ and trilinear operator $b(u,v,w)$ in have the properties
\begin{equation}
\left\{\aligned
	& \|B(u,u)\|^2_{V'} \leq c_0|u|\|u\|, \quad \forall \, u \in V, \\
	& b(u,v,v) = 0, \quad \forall \, u \in V, \; \forall \, v \in (H^1_0(\Omega))^3,\\
	& b(u,v,w) \leq C|u|^{1/4}\|u\|^{3/4}\|v\||w|^{1/4}\|w\|^{3/4}, \quad \forall \, u, v, w \in V, \\
	& b(u,v,w) = -b(u,w,v), \quad \forall \, u,v,w \in V.
\endaligned\right.\label{q7}
\end{equation}
\end{Lemma}

\begin{Lemma}\label{le4.8} (1) Assume $f(x)\in V'$ and $u_{0}\in H$ in \eqref{bmy1-1}, then the semigroup $\{S(t)\}$ has
	a bounded absorbing ball $B_0=\left\{u\in H:|u|\leq  \rho\right\}$ in $H$.
	
	(2) If $f(x)\in H$ and $u_{0}\in H$ in \eqref{bmy1-1}, then the semigroup $\{S(t)\}$ has
	a bounded absorbing ball $\hat{B}_0=\{u\in D(A^{1/2}):\|u\|_{D(A^{{1}/{2}})}\leq  \hat{\rho}\}$ in  $D(A^{{1}/{2}})$.
\end{Lemma}

\noindent{\bf Proof.}
	(1) Taking inner product of \eqref{bmy1-1} with $u$ and integrating by parts over $\Omega$, we derive
	\begin{align} \frac{1}{2}\frac{d}{d t}|u|^2+\nu\| u\|^2+\nu_0\|u\|^4
	= &\langle f(x),u\rangle  \nonumber\\
	\leq& \frac{8}{\nu}\|f(x)\|^2_{V'}+\frac{\nu}{2}\|u\|^2,\label{25}
	\end{align}
	since $(B(u,u),u)=0$.
	
	From the Poincar\'{e} inequality and Gronwall's inequality, \eqref{25} yields
	\begin{eqnarray}
	|u(t,x)|^2\leq e^{-\nu\lambda_1 (t-\tau)}|u_0|^2+\frac{8}{\nu^2\lambda_1}\|f(x)\|^2_{V'}.\label{b11}
	\end{eqnarray}
	Choosing a time $T_0=\tau+\frac{1}{\nu\lambda_1}\ln\Big(\frac{\nu^2_0\lambda_1|u_0|^2}{8\|f(x)\|_{V'}^2}\Big)+1$ such that $e^{-\nu\lambda_1 (t-\tau)}|u_0|^2\leq \frac{1}{\nu^2\lambda_1}\|f(x)\|^2_{V'}$ if $t\geq T_0$. Defining $\rho^2=\frac{16}{\nu\lambda_1}\|f(x)\|^2_{V'}$, we conclude that $B_0=\Big\{u:|u|\leq \rho\Big\}$ is the bounded absorbing ball in $H$.
	
	Moreover, integrating from $t$ to $t+1$ with respect to time variable for \eqref{b11}, we can derive the estimate
	\begin{eqnarray}
	\nu \int^{t+1}_t\|u(s)\|^2ds+2\nu_0\int^{t+1}_t\|u(s)\|^4ds\leq \frac{16}{\nu}\|f\|^2_{V'}+|u(t)|^2.\label{xin-1}
	\end{eqnarray}
	
(2) Taking inner product of \eqref{bmy1-1} with $A u$ and integrating by parts over $\Omega$, we derive
	\begin{equation} \frac{1}{2}\frac{d}{d t}|A^{1/2}u|^2+\nu| A u|^2+\nu_0\|u\|^2|A u|^2\leq b(u,u,Au)+\langle f(x),A u\rangle,\label{25-1}\end{equation}
	From the Cauchy-Schwarz inequality, the generalized Poincar\'{e} inequality and the property of trilinear operator, we have
$$
\langle f(x),A u\rangle \leq |f(x)||A u| \leq \frac{\nu}{2}|A u|^2+\frac{C}{\nu}|f(x)|^2 ,
$$
\begin{align}
|b(u,u, A u)|
& \leq  \int_{\Omega}|u||\nabla u||A u|d x\leq \|u\|_{L^4}\|\nabla u\|_{L^4}|Au|\nonumber\\
& \leq  C|u|^{\frac{1}{4}}|\nabla u|^{\frac{3}{4}}|\nabla u|^{\frac{1}{4}}
|Au|^{\frac{3}{4}}|A u|=C|u|^{\frac{1}{4}}|Au|^{\frac{3}{4}}|\nabla u||Au|\nonumber\\
& \leq \nu_0\|u\|^2|Au|^2+\frac{\nu}{2}|Au|^2+C|u|^2. \label{yang-4}
\end{align}
	Combining \eqref{25-1}-\eqref{yang-4}, integrating the time variable from $s$ to $t$, here $t-1\leq s\leq t$, it yields
	\begin{eqnarray}
	|A^{1/2}u(t)|^2\leq |A^{1/2}u(s)|^2+\frac{C}{\nu}|f(x)|^2+C\|u\|^2_{L^{\infty}(\tau,T;H)},\label{xin-2}
	\end{eqnarray}
	then integrating with respect to $s$ from $t-1$ to $t$ for \eqref{xin-2}, using \eqref{xin-1}, we derive that
\begin{eqnarray*}
	|A^{1/2}u(t)|^2
	&\leq& \int^t_{t-1}|A^{1/2}u(s)|^2ds+\frac{C}{\nu}|f(x)|^2+C\rho^2\\
	&\leq& \frac{16}{\nu^2}\|f\|^2_{V'}+\frac{1}{\nu}|u(t)|^2+\frac{C}{\nu}|f(x)|^2+C\rho^2\\
	&\leq&\frac{16}{\nu^2}\|f\|^2_{V'}+\frac{C}{\nu}|f(x)|^2+(\frac{1}{\nu}+C)\rho^2\\
	&=& \hat{\rho}^2,
\end{eqnarray*}
	we conclude that $\hat{B}_0=\Big\{u:|u|\leq \hat{\rho}\Big\}$ is the bounded absorbing ball in $D(A^{\frac{1}{2}})$ for $t\geq T_0$. \qed

\begin{Lemma}\label{le4.11} For any $f(x)\in V'$ and $u_0\in H$, the
	semigroup $\{S(t)\}$ generated
	by the system \eqref{yang-2} is asymptotically compact in $H$.
\end{Lemma}

\noindent{\bf Proof.} Since the embedding $D(A^{\frac{1}{2}})\hookrightarrow\hookrightarrow H$ is compact, we can deduce the asymptotic compactness for the semigroup. Combining with Lemma \ref{le4.8}, we complete the proof. \qed \\

\begin{Lemma}\label{le4.8-1}  If $f(x)\in H$ and $u_{0}\in D(A^{\frac{1}{4}})$ in \eqref{bmy1-1}, then the semigroup $\{S(t)\}$ has
	a bounded absorbing ball $\bar{B}_0=\{u\in D(A^{\frac{1}{4}}):\|u\|_{D(A^{\frac{1}{4}})}\leq  \bar{\rho}\}$ in  $D(A^{\frac{1}{4}})$.
\end{Lemma}

\noindent{\bf Proof.}
Since the domain $\Omega$ is bounded and $u_{0}\in D(A^{\frac{1}{4}})$, by Lemma \ref{le4.8}.  we can easily deduce the existence of bounded absorbing ball in $D(A^{\frac{1}{4}})$ with radius $\bar{\rho}$. \qed

\begin{Lemma}\label{le4.12}
	Let $f(x)\in H$ and $u_0\in D(A^{\frac{1}{4}})$. Then the semigroup $\{S(t)\}$ generated
	by the problem \eqref{yang-2} satisfies (MWZ). In particular it is asymptotically smooth in $D(A^{\frac{1}{4}})$.
\end{Lemma}

\noindent{\bf Proof.} {\it Step 1:}
From Lemmas \ref{le4.8} and \ref{le4.8-1}, we see that $\hat{B}_0$ be the bounded absorbing set, then there exists a forward time $t_{\hat{B}_0}$ such that
	$\|S(t)u_0\|_{D(A^{\frac{1}{4}})}^2\leq \bar{\rho}$.
	
		From the definition in Section 2, that $D(A^{\frac{1}{4}})=\{\omega'_1,\omega'_2,\cdots, \omega'_m, \cdots\}$ is a Hilbert space with orthonormal decomposition
		$$
		D(A^{\frac{1}{4}})=D_1(A^{\frac{1}{4}})\oplus D_2(A^{\frac{1}{4}}),
		$$
		where $D_1(A^{\frac{1}{4}})=span\{\omega'_1,\omega'_2,\cdots, \omega'_m\}$,  $D_2(A^{\frac{1}{4}})=span\{\omega'_{m+1},\omega'_{m+2},\cdots \}$.  Let $P$ be an orthonormal projector from $D(A^{\frac{1}{4}})$ to $D_1(A^{\frac{1}{4}})$.
 Hence the solution $u$
	has the decomposition
	\begin{equation}
	u=Pu+(I-P)u:=u_1+u_2,
	\end{equation} for $u_1\in D_1(A^{\frac{1}{4}})$, $u_2\in D_2(A^{\frac{1}{4}})$ with the initial data $A^{\frac{1}{4}}u_1(\tau)=PA^{\frac{1}{4}}u_0$ and $A^{\frac{1}{4}}u_2(\tau)=(I-P)A^{\frac{1}{4}}u_0$ (which also are bounded in $D(A^{\frac{1}{4}})$) respectively.
	
	{\it Step 2:} Since $u_1$ is the orthonormal projection of $u$, from the existence of absorbing ball $\hat{B}_0$ for the semigroup $S(t)$, we derive that $u$ is bounded in $D(A^{\frac{1}{4}})$, and hence $u_1$ is bounded in $D(A^{\frac{1}{4}})$, i.e., $|A^{\frac{1}{4}}u_1|^2\leq \bar{\rho}$.

	{\it Step 3:} The next objective is to obtain the $D(A^{\frac{1}{4}})$-norm of $u_2$ is small enough as $m\rightarrow +\infty$.
	
	By taking inner product of \eqref{yang-2} with $A^{1/2} u_2$ and noting that $(A^{\frac{1}{4}}u_1, A^{\frac{1}{4}}u_2)=0$ in $D(A^{\frac{1}{4}})$ and $(A^{1/2}u_1, A^{1/2} u_2)=0$ by projector $P$,
	we obtain
	\begin{align}
	\frac{1}{2}\frac{d}{dt}|A^{\frac{1}{4}}u_2|^2
	& + \nu|A^{3/4} u_2|^2 +\nu_0\|u\|^2|A^{3/4} u_2|^2\leq  |(B(u,u),A^{1/2} u_2)|+|\langle P f,A^{1/2} u_2\rangle|. \label{g2-1}
	\end{align}
By the property of $b(\cdot,\cdot,\cdot)$ in \eqref{q7}, using the $\varepsilon$-Young inequality, we have
	\begin{align*}
	 |(B(u,u),A^{1/2} u_2)|
	 & \leq \int_{\Omega}|u||A^{1/2} u||A^{1/2} u_2|d x\leq C|u|_{L^4}\|u\|_{L^4}|A^{1/2} u_2|\\
	 &\leq C|u|^{\frac{1}{4}}\|u\|^{\frac{3}{4}}\|u\|^{\frac{1}{4}}|Au|^{\frac{3}{4}}  \frac{1}{\lambda^{1/4}_{m+1}}|A^{3/4}u_2| \\
	 & \leq \nu_0\|u\|^2|A^{3/4}u_2|^2+\frac{C_1}{\nu_0\lambda^{1/4}_{m+1}}\Big(|u|^{2}+C_2|Au|^{2}\Big).
	 \end{align*}
	Using the definition in \eqref{2} and Young's inequality, we have
	\begin{equation}
	|\langle P(f-f^{\varepsilon}),A^{1/2} u_2\rangle|\leq \frac{C}{\nu\lambda_{m+1}^{1/4}}|f|^2+\frac{\nu}{2}|A^{3/4} u_2|^2.\label{g2-6}
	\end{equation}
	Combining \eqref{g2-1}--\eqref{g2-6}, we conclude
$$
\frac{d}{d t}|A^{\frac{1}{4}}u_2|^2+\nu|A^{3/4} u_2|^2\leq \frac{C}{\nu\lambda_{m+1}^{1/4}}|f|^2+\frac{C_1}{\nu_0\lambda^{1/4}_{m+1}}\Big(|u|^{2}+C_2|Au|^{2}\Big).
$$
Using the Poincar\'{e} inequality and noting $A^{\alpha}u_2=\displaystyle{\sum^{\infty}_{i=m+1}}\lambda^{\alpha}_i (u_2,\omega_i)\omega_i$, we conclude
\begin{equation}\frac{d}{d t}|A^{\frac{1}{4}}u_2|^2+\nu\lambda_{m+1}|A^{\frac{1}{4}}u_2|^2
	 \leq \frac{C}{\nu\lambda_{m+1}^{1/4}}|f|^2+\frac{C_1}{\nu_0\lambda^{1/4}_{m+1}}\Big(|u|^{2}+C_2|Au|^{2}\Big),\label{gg10-1}
\end{equation}
by applying the Gronwall inequality in $[\tau,t]$ to \eqref{gg10-1} and the H\"{o}lder inequality, we deduce that
	\begin{align}
	|A^{\frac{1}{4}}u_2(t)|^2
\leq & \, |A^{\frac{1}{4}}u_2(\tau)|^2e^{-\nu\lambda_{m+1}(t-\tau)}+\frac{C}{\nu\lambda_{m+1}^{1/4}}|f|^2\nonumber \\
& \, + \frac{C_1}{\nu_0\lambda^{1/4}_{m+1}}\int^t_{\tau}(C_1|u(s)|^{2}+C_2|Au|^{2})e^{-\nu\lambda_{m+1}(t-s)}d s\nonumber\\
\leq & \, e^{-\nu\lambda_{m+1}(t-\tau)}\hat{\rho}^2+\frac{C|f|^2}{\nu\lambda_{m+1}^{1/4}}\nonumber\\
& \, +\frac{C}{\nu_0\lambda^{1/4}_{m+1}}\|u\|^{2}_{L^{\infty}(\tau,T;H)}(1-e^{-\nu\lambda_{m+1}(t-\tau)})\nonumber\\
& \, +\frac{C}{\nu_0\lambda^{1/4}_{m+1}}\|Au\|^{2}_{L^2(\tau,t;D(A))}(1-e^{-\nu\lambda_{m+1}(t-\tau)})\nonumber\\
= & \, I_1+I_2+I_3+I_4.\label{g7-1}
	\end{align}
	By the boundedness of absorbing set, noting $\displaystyle{\lim_{m\rightarrow\infty}}\lambda_{m+1}=+\infty$ and the existence of global solution, then for $m$ large enough, it follows that
\begin{align*}
I_1 & =|A^{\frac{1}{4}}u_2(\tau)|^2e^{-\nu\lambda_{m+1}(t-\tau)}  \leq  \frac{\varepsilon}{4}, \nonumber \\
I_2 & = \frac{C|f|^2}{\nu\lambda_{m+1}^{1/4}}  < \frac{\varepsilon}{4}, \nonumber \\
I_3 & \leq \frac{C}{\nu_0\lambda^{1/4}_{m+1}}\|u\|^{2}_{L^{\infty}(\tau,T;H)}  < \frac{\varepsilon}{4}, \nonumber \\
I_4 & \leq \frac{C}{\nu_0\lambda^{1/4}_{m+1}}\|Au\|^{2}_{L^2(\tau,t;D(A))}  < \frac{\varepsilon}{4}.
\end{align*}
Combining  \eqref{g7-1} with above estimates we conclude that
$$
\|(I-P)S(t)u_{\tau}\|^2_{D(A^{\frac{1}{4}})}=|A^{\frac{1}{4}}u_2|^2< \varepsilon ,
$$
which implies the (MWZ) condition. The result follows by using Lemma \ref{le4.8}. \qed

\subsection{Proof of Theorem \ref{th5.1}} \label{se3.2}
 In this section, we prove the existence of pullback attractors in $H$ and $D(A^{\frac{1}{4}})$ for \eqref{bmy1-1}.

\begin{Lemma}\label{le5.2}
Assume the external force $f(x,t)\in L^2_{loc}(\mathbb{R};V')$ satisfies \eqref{f7a} and fix a parameter $\mu\in (0,\mu_0]$ $(\mu_0=\nu\lambda_1)$.
Then the solution $u$ of problem \eqref{bmy1-1}, with initial data $u_0\in H$, satisfies for any $\tau\leq t$,
$$
	|u|^2\leq |u_0|^2e^{-\mu(t-\tau)}+\frac{Ce^{-\mu t}}{\nu}\int^t_{-\infty}e^{\mu s}\|f(s)\|^2_{V'}d s.
$$
\end{Lemma}

\noindent{\bf Proof.} Multiplying \eqref{bmy1-1} with $u$ and integrating over $\Omega$, choosing an appropriate parameter $\mu$ such that $e^{-\mu(t-\tau)}\geq e^{-\mu_0(t-\tau)}$ and $$\frac{e^{-\mu t}}{\nu}\int^t_{-\infty}e^{\mu s}\|f(s)\|^2_{V'}d s \geq \frac{e^{-\mu_0 t}}{\nu}\int^t_{-\infty}e^{\mu_0 s}\|f(s)\|^2_{V'}d s,$$ we can derive the result easily. \qed

\begin{Lemma}[Pullback $\mathcal{D}$-absorbing in $H$] \label{le5.3}
Assume $f(t,x) \in L^2_{loc}(\mathbb{R};V')$ satisfies \eqref{f7a} and
let $\hat{B}_0=\{B_0(t)\}_{t\in\mathbb{R}}$ be a family of balls,
where
%$B_0(t)=\overline{B(0,\rho_0(t))}$ is a ball at center $0$ and radius $\rho_0(t)$ that satisfies
$$
\rho^2_0(t)=1+\frac{4}{\nu\lambda_1}\int^t_{-\infty}e^{-\mu(t-s) }\|f(s)\|^2_{V'}ds.
$$
Then for any $0<\varepsilon<\frac{1}{2}$ small enough, there exists a pullback time $\tau(t,\varepsilon)$, such that for any $\tau< \tau(t,\varepsilon)\leq t$, $\hat{B}_0(t)$ is a family of pullback $\mathcal{D}$-absorbing sets for the continuous  process $S(t,\tau): H \to H$.
\end{Lemma}

\noindent{\bf Proof.} Noting that
$$
|S(t,\tau)u_0|^2 \leq |u_0|^2e^{-\mu(t-\tau)}+\frac{Ce^{-\mu t}}{\nu}\int^t_{-\infty}e^{\mu s}\|f(s)\|^2_{V'}ds,
$$
there exists a pullback time $\tau(t,\varepsilon)$, such that for $\tau<\tau(t,\varepsilon)$, if follows
$e^{-\mu t}|\rho_{\hat{D}}(\tau)|^2e^{\mu\tau}\leq \varepsilon$.
Hence, we have
$$
|S(t,\tau)u_0|^2\leq \varepsilon+\rho^2_0(t)-\frac{1}{2}\leq \rho^2_0(t),
$$
which implies that $S(t,\tau)D(\tau)\subset B_0(t)$. Then $\hat{B}_0(t)$ is a family of pullback $\mathcal{D}$-absorbing balls. \qed

\medskip

The pullback $\mathcal{D}$-asymptotic compactness for the process in $H$ is obtained by verifying (MWZ) condition.

\begin{Lemma}[Pullback $\mathcal{D}$-asymptotic compactness in $H$] \label{le5.5}
Assume $f(x,t)$ belongs to $L^2_{loc}(\mathbb{R};V')$ which satisfies \eqref{f7a}.
Then the process $S(t,\tau)$ is pullback $\mathcal{D}$-asymptotically compact in $H$ for the system \eqref{bmy1-1}.
\end{Lemma}

\noindent{\bf Proof.} The pullback $\mathcal{D}$-asymptotic compactness follows from the (MWZ) condition as in Lemma \ref{le4.12}. Here we omit the details. \qed

\begin{Lemma}\label{le5.4}
	Assume $f(x,t)\in L^2_{loc}(\mathbb{R};H)$ satisfies \eqref{f7a-1}, and fix a parameter $\mu\in (0,\mu_0]$ $(\mu_0=\nu\lambda_1)$.
	Then the solution $u$ of problem \eqref{bmy1-1}, $u_{0}\in D(A^{\frac{1}{2}})$, satisfies for any $\tau\leq t$,
	\begin{align}
	|A^{\frac{1}{4}}u(t)|^2
	 \leq & \; |A^{\frac{1}{4}}u_0|^2e^{-\mu(t-\tau)}+\frac{Ce^{-\mu t}}{\nu}\int^t_{-\infty}e^{\mu s}|f(s)|^2d s\nonumber\\  &+\frac{C}{\nu_0}\Big(\|u\|^2_{L^{\infty}(\tau,T;H)}+\|u\|^2_{L^2(\tau,T;D(A))}\Big).\label{ya-7}
	\end{align}
\end{Lemma}

\noindent{\bf Proof.} Multiplying \eqref{bmy1-1} with $A^{1/2} u$ and integrating over $\Omega$,  we derive that
$$
\frac{1}{2}\frac{d}{dt}|A^{\frac{1}{4}}u|^2 + \nu|A^{3/4} u|^2+\nu_0\|u\|^2|A^{3/4} u|^2 \leq  |b(u,u,A^{1/2} u)| + | \langle P f(t,x),A^{1/2} u \rangle|. \label{ya2-1}
$$
Then using the similar technique in \ref{le4.8}, from the Poincar\'{e} inequality and Gronwall's inequality,
	choosing an appropriate parameter $\mu$ such that $e^{-\mu(t-\tau)}\geq e^{-\mu_0(t-\tau)}$ and $$\frac{e^{-\mu t}}{\nu}\int^t_{-\infty}e^{\mu s}|f(s)|^2d s \geq \frac{e^{-\mu_0 t}}{\nu}\int^t_{-\infty}e^{\mu_0 s}|f(s)|^2d s,$$
	it follows that \eqref{ya-7} holds. The proof is completed. \qed

\begin{Lemma}[Pullback $\mathcal{D}$-absorbing in $D(A^{\frac{1}{4}})$] \label{le5.4.1}
	Assume the external force $f(x,t)\in L^2_{loc}(\mathbb{R};H)$ satisfies \eqref{f7a-1}, and let $\hat{B}'_0=\{B'_0(t)\}_{t\in\mathbb{R}}$ be a family of balls, where $B'_0(t)=\overline{B'(0,\rho'_0(t))}$ is a ball with
$$
(\rho'_0(t))^2=1+\frac{C}{\nu_0}\Big(\|u\|^2_{L^{\infty}(\tau,T;H)}+\|u\|^2_{L^2(\tau,T;D(A))}\Big)
+\frac{C}{\nu}\int^t_{-\infty}e^{-\mu(t-s)} |f(s)|^2 ds.
$$
Then for any $0<\varepsilon<\frac{1}{2}$ small enough, there exists a pullback time $\tau'(t,\varepsilon)$, such that for any $\tau< \tau'(t,\varepsilon)\leq t$, $\hat{B}'_0(t)$ is a family of pullback $\mathcal{D}$-absorbing sets for the continuous  process $S(t,\tau):D(A^{\frac{1}{4}}) \to D(A^{\frac{1}{4}})$.
\end{Lemma}

\noindent{\bf Proof.} Noting Lemma \ref{le5.4},
	there exists a pullback time $\tau'(t,\varepsilon)$, such that for any $\tau<\tau'(t,\varepsilon)$, if it follows that
	$e^{-\mu t}|\rho'_{\hat{D}'}(\tau)|^2e^{\mu\tau}\leq \varepsilon$.
Hence, we have
$$
\|S(t,\tau)u_0\|^2_{D(A^{\frac{1}{4}})}\leq \varepsilon+(\rho'_0(t))^2-\frac{1}{2}\leq (\rho'_0(t))^2,
$$
which implies that $S(t,\tau)D'(\tau)\subset B'_0(t)$. That is, $\hat{B}'_0(t)$ is a family of pullback $\mathcal{D}$-absorbing balls in $D(A^{\frac{1}{4}})$. \qed

\begin{Lemma}[Pullback $\mathcal{D}$-asymptotic compactness in $D(A^{\frac{1}{4}})$] \label{le5.6}	
	Assume the external force $f(t,x)\in L^2_{loc}(\mathbb{R};H)$ satisfies \eqref{f7a-1}.
	Then the process $S(t,\tau): D(A^{\frac{1}{4}}) \to D(A^{\frac{1}{4}})$ satisfies condition (MWZ).
	In particular the process is pullback $\mathcal{D}$-asymptotically compact in $D(A^{\frac{1}{4}})$ for problem \eqref{bmy1-1}.
\end{Lemma}

\noindent{\bf Proof.}
	{\it Step 1:} Let $\mathcal{B}'=\{B'_0(t)\}_{t\in\mathbb{R}}$ be the pullback $\mathcal{D}$-absorbing family given in Lemma \ref{le5.4.1}. Then there exists a pullback time $\tau_{t,\varepsilon_1}$ such that
	$|A^{1/4}u|^2=|S(t,\tau)u_0|^2\leq \rho'_{0}(t)$.
	
The decomposition is similar to Lemma \ref{le4.12}, hence we can write solution $u$ as
$$
u=Pu+(I-P)u:=u_1+u_2,
$$
for $u_1\in D_1(A^{\frac{1}{4}})$, $u_2\in D_2(A^{\frac{1}{4}})$ with the initial data which also are bounded in $D(A^{\frac{1}{4}})$.
	From the existence of global solution and pullback $\mathcal{D}$-absorbing family of set, we know $\|A^{1/4}u_1\|^2\leq \rho'_{0}(t)$, and next we only need to prove the $D(A^{\frac{1}{4}})$-norm of $u_2$ is small enough.
	
	{\it Step 2:} Taking inner product of \eqref{bmy1-1} with $A^{1/2} u_2$, using the same technique in Lemma \ref{le4.12},
we conclude
\begin{equation} \label{gg10-1-1}
 \frac{d}{dt}|A^{\frac{1}{4}}u_2|^2+\nu\lambda_{m+1}|A^{\frac{1}{4}}u_2|^2
 \leq \frac{C}{\nu\lambda_{m+1}^{1/4}}|f|^2+\frac{C_1}{\nu_0\lambda^{1/4}_{m+1}}\Big(|u|^{2}+C_2|Au|^{2}\Big).
\end{equation}
By applying the Gronwall inequality in $[\tau,t]$ to \eqref{gg10-1-1}, we deduce that
\begin{align}
|A^{\frac{1}{4}}u_2(t)|^2
\leq & \; e^{-\nu\lambda_{m+1}(t-\tau)}\hat{\rho}^2
+ \frac{C}{\nu\lambda_{m+1}^{1/4}}\int^t_{-\infty}e^{-\mu(t-s)}|f(s)|^2ds \nonumber\\
& \, +\frac{C}{\nu_0\lambda^{1/4}_{m+1}}\|u\|^{2}_{L^{\infty}(\tau,T;H)}
+\frac{C}{\nu_0\lambda^{1/4}_{m+1}}\|Au\|^{2}_{L^2(\tau,t;D(A))}\nonumber\\
& = \, J_1+J_2+J_3.\label{g7-1-1}
\end{align}
From the existence of bounded family of pullback $\mathcal{D}$-absorbing sets, noting that $\displaystyle{\lim_{m\rightarrow\infty}}\lambda_{m+1}=+\infty$, $\tau\rightarrow -\infty$, $f(s)\in L^2_{loc}(\mathbb{R};H)$ and the existence of global solution, then for $m$ large enough, there exists a pullback time $\tau_0$, such that for all $\tau\leq \tau_0$,
\begin{align*}
	J_1=|A^{\frac{1}{4}}u_2(\tau)|^2e^{-\nu\lambda_{m+1}(t-\tau)}\leq (\rho'_{0}(t))^2 e^{\nu\lambda_{m+1}\tau} & < \frac{\varepsilon}{3},\\
	 J_2=\frac{C}{\nu\lambda_{m+1}^{1/4}}\int^t_{-\infty}e^{-\mu(t-s)}|f(s)|^2ds & < \frac{\varepsilon}{3},\\
	J_3+J_4 & < \frac{\varepsilon}{3}.
\end{align*}
Combining \eqref{g7-1-1} with above estimates we conclude
$$
 \|(I-P)S(t,\tau)u_{\tau}\|^2_{D(A^{\frac{1}{4}})}< \varepsilon,
$$
which implies the process satisfies (MWZ) condition.
Hence $\{S(t,\tau)\}$ is pullback $\mathcal{D}$-asymptotically compact in $D(A^{\frac{1}{4}})$.  \qed

\medskip

\noindent{\bf Proof of Theorem \ref{th5.1}.} We proved the existence of families of pullback $\mathcal{D}$-absorbing sets in $H$ and $D(A^{\frac{1}{4}})$ by Lemmas \ref{le5.3} and \ref{le5.4.1}. We also proved the pullback asymptotic compactness for  continuous processes $S(t,\tau)$ in $H$ and $D(A^{\frac{1}{4}})$, in Lemmas \ref{le5.5} and \ref{le5.6} respectively.
Then our result follows from Theorem \ref{thm-EPA}. \qed

\subsection{Proof of Theorem \ref{th5-f}} \label{sec4.3}

Let us consider the Cauchy problem, variation of \eqref{yang-2},
\begin{equation}
\left\{\aligned &\frac{d U}{d t}+\nu AU+\mathbb{A}U+B(u,U)+B(U,u)=0,\\
&U(\tau)=\xi\in H.\endaligned\right.\label{5-5-1}
\end{equation}
Taking inner product of \eqref{5-5-1} by $U$, we can prove the existence of unique solution
$U\in L^2(\tau,T;V)\cap L^4(\tau,T;V)\cap C(\tau,+\infty;H)$ for any $T>\tau$ and $\tau\in \mathbb{R}$.

\begin{Lemma}\label{le4.9.1}
	The solution of \eqref{5-5-1} generates a bounded linear compact operator
$$
\Lambda(t,s;u_0)\xi=U(t): H\rightarrow H,
$$
that satisfies
	\begin{equation}
	\sup_{\hat{u}_0, u_0\in \mathcal{A}(s)}\sup_{\|\hat{u}_0-u_0\|\leq \varepsilon}\frac{\|S(t,s)\hat{u}_0-S(t,s)u_0-\Lambda(t,s;\hat{u}_0)(\hat{u}_0-u_0)\|}{\|\hat{u}_0-u_0\|}\rightarrow 0, \label{5-5-9}
	\end{equation}
as $\varepsilon\rightarrow 0$. This means the process is uniformly differentiable for $t\geq s$.
\end{Lemma}

\noindent{\bf Proof.} (1) Assume $w=u-\hat{u}$, $u(t)$ and $\hat{u}(t)$ be two solutions of
$$
\frac{du}{dt}+\nu Au+\mathbb{A}u+B(u,u)=Pf(t)
$$
with different initial data $u(s)=u_0$ and $\hat{u}(s)=\hat{u}_0$.
	Denoting $U(t)$ be a solution of problem \eqref{5-5-1} with initial data $U(s)=u_0-\hat{u}_0$, we can verify that $\theta=u-\hat{u}-U$ satisfies the Cauchy problem
\begin{equation}
	\left\{\aligned &\frac{d\theta}{d t}+\nu A\theta+\mathbb{A}\theta+B(u,\theta)+B(\theta,u)-B(w,w)=0,\\
	&\theta(s)=0.\endaligned\right.\label{5-5-2}
\end{equation}
Taking inner product of \eqref{5-5-2} by $\theta$, using the property of operator $\mathbb{A}$, we can derive
	\eqref{5-5-9} easily.
	
	(2) Taking inner product of \eqref{5-5-1} with $U$ and $A U$, integrating over $\Omega$, by the same technique in the proof of regular pullback absorbing set, we can derive the uniform estimate in more regular space $D(A^{\frac{1}{2}})$. Since $D(A^{\frac{1}{2}})$ is compact in $H$, then we can prove that the operator $\Lambda(t,s;u_0)$ is compact. \qed

\medskip

Next, we shall use trace formula \cite[Lemmas 4.19 and 4.20]{clr} to prove the boundedness of fractal dimension of pullback attractors in $H$.

\medskip

\noindent{\bf Proof of Theorem \ref{th5-f}.} From the definition of the family of pullback attractors, then for a fixed $\tau^*$, $\displaystyle{\bigcup_{\tau\leq \tau^*}}{\mathcal{A}}(t)$ is precompact in $H$.
For each $t\geq \tau$, $\tau\leq \tau^*$ and $u_0\in H$, the linear operator is described as $\Lambda(t,s;u_0)\cdot \xi=U(t)$, where $U(t)$ is the solution of \eqref{5-5-1}.
Denoting $F(S(t,\tau)u_0,t)=-\nu A-\mathbb{A}-B(u,\cdot)-B(\cdot,u)$, then from Lemma \ref{le4.9.1}, we see that $F(\cdot,t)$ is Gateaux differential in $V$ at $S(t,\tau)u_0$ which satisfies
$$
F'(S(t,\tau)u_0)U=-\nu A U-\mathbb{A} U-B(S(t,\tau)u_0,U)-B(U,S(t,\tau)u_0),
$$
this implies $F'(S(t,\tau)u_0,t)\in \mathcal{L}(V,V')$ is a continuous linear operator satisfying the problem
\begin{equation}
\left\{\aligned &\frac{d U}{d t}=F'(S(t,\tau)u_0,t)U,\ \ u_0\in H,\\
&U(\tau,x)=\xi\endaligned\right.\label{yan-11}
\end{equation}
which possesses a unique solution $U(t)=U(t,\tau;u_0,\xi)\in L^2(\tau,T;V)\cap C(\tau,T;H)$.

For each $\xi_1,\ \xi_2,\ \cdots, \xi_n\in H$, we denote $U_i(t)=\Lambda(t,s;u_0)\cdot \xi_i$ which implies $U_1(s)=U_1(s,\tau;u_0,\xi_1)$, $ U_2(s)=U_2(s,\tau;u_0,\xi_2)$, $\cdots$, $U_n(s)=U_n(s,\tau;u_0,\xi_n)$ be the solution of problem \eqref{yan-11} with different initial data $U_i(\tau)=\xi_i$ $(i=1,2,\cdots,n)$ respectively, $Q_n(s)$ denote the projection from $H$ to the space $span\{U_1(s),\ U_2(s),\ \cdots,\ U_n(s)\}$, then by Lemma 4.19 in \cite{clr},
it yields
\begin{align*}
\|U_1(t) & \wedge U_2(t) \wedge \cdots \wedge U_n(t)\|_{\wedge^n(H)} \\
& = \|\xi_1\wedge \xi_2\wedge\cdots \wedge \xi_n\|_{\wedge^n(H)}\exp\Big(\int^t_{\tau}\mbox{Tr}_n(F'(S(s,\tau)u_0,s)\circ Q_n(s)d s)\Big).
\end{align*}
Let $\{e_1(s),\ e_2(s),\ \cdots,\ e_n(s)\}$ be an orthonormal basis for $span\{U_1(s)$, $U_2(s)$, $\cdots$, $U_n(s)\}$,
then
$$
\mbox{Tr}_n(F'(S(s,\tau)u_0,s)=\sup_{\xi_i\in H, |\xi_i|\leq 1, i\leq n}\Big(\sum^n_{i=1} \langle F'(S(s,\tau)u_0,s)e_i,e_i\rangle\Big).\label{yan-13}
$$
Since $U_i(s)\in L^2(\tau,T;V)$, then $U_i(s)\in V$ for a.e. $s\geq \tau$, hence $e_i(s)\in V$ for a.e. $s\geq \tau$ and $i=1,2,\cdots,n$.

Noting that $b(S(t,\tau)u_0,e_i(s),e_i(s))=0$, we derive
\begin{align}
\mbox{Tr}_n & (F'(S(s,\tau)u_0,s)\circ Q_n(s)\nonumber\\
& = \sum^n_{i=1}\langle F'(S(s,\tau)u_0,s)e_i(s),e_i(s)\rangle\nonumber\\
& = \sum^n_{i=1}\Big(-\nu\|e_i(s)\|^2-b(S(s,\tau)u_0,e_i(s),e_i(s))-b(e_i(s),S(s,\tau)u_0,e_i(s))\Big)\nonumber\\
& \leq  -\nu\sum^n_{i=1}\|e_i(s)\|^2-\nu_0\sum^n_{i=1}\|e_i(s)\|^4+\sum^n_{i=1}|b(e_i(s),S(s,\tau)u_0,e_i(s))|. \label{yan-14}
\end{align}
For the second term in \eqref{yan-14},  by the Lieb-Thirring inequality in 3D case ($p=2, \ n=3$):
$$
\Big(\int_{\Omega}\Big(\displaystyle{\sum^n_{i=1}}|e_i(s)|^2\Big)^{\frac{p}{p-1}}d x\Big)^{\frac{2(p-1)}{n}}\leq C_1\displaystyle{\sum^n_{i=1}}\int_{\Omega}|\nabla e_i(s)|^2d s,
$$
where $\frac{n}{2}<p\leq 1+\frac{n}{2}$, we could proceed using the bound
\begin{align*}
\sum^n_{i=1}|b(e_i(s),S(s,\tau)u_0,e_i(s))|
& \leq C\int_{\Omega}\Big(\sum^n_{i=1}|e_i(s)||\nabla S(s,\tau)u_0||e_i(s)|\Big)d x \\
& \leq \|S(s,\tau)u_0\| \Big[\int_{\Omega}\Big(\sum^n_{i=1}|e_i(s,x)|^2\Big)^{2}d s\Big]^{\frac{2}{3}\times\frac{3}{4}} \\
& \leq \frac{C}{\nu}\|S(s,\tau)u_0\|^2+\frac{\nu}{2}\Big[\sum^n_{i=1}\|e_i(s)\|^2\Big]^{\frac{3}{4}}.
\end{align*}
Using the variational principle and $\displaystyle{\sum^n_{i=1}\lambda_i}\geq \frac{\pi n^2}{|\Omega|}$ from \cite{ily}, taking the average, we obtain
\begin{align*}
 {\rm Tr}_n (F'(S & (s,\tau)u_0,s) \circ Q_n(s)\\
& \leq  -\nu\sum^n_{i=1}\|e_i(s)\|^2-\nu_0\sum^n_{i=1}\|e_i(s)\|^4+\frac{\nu}{2}\Big[\sum^n_{i=1}\|e_i(s)\|^2\Big]^{\frac{3}{4}}+\frac{C}{\nu} \|S(s,\tau)u_0\|^2\\
& \leq \frac{2\nu}{27}-\frac{\nu}{2}\sum^n_{i=1}\lambda_i-\nu_0\sum^n_{i=1}\lambda_1^2+\frac{C}{\nu}\|S(s,\tau)u_0\|^2\\
& \leq \frac{2\nu}{27}-\frac{\pi\nu n^2}{2|\Omega|}-\frac{\pi^2\nu_0 n^4}{|\Omega|^2}+\frac{C}{\nu} \|S(s,\tau)u_0\|^2.
\end{align*}
Defining
\begin{align*}
&q_n=\sup_{t\in\mathbb{R}}\sup_{u_0\in\mathcal{A}(t)}\Big(\frac{1}{T}\int^t_{t-T}Tr_n(F'(S(s,\tau)u_0,s)\circ Q_n(s)d s\Big),\\
&\hat{q}_n=\limsup_{T\rightarrow +\infty}q_n,
\end{align*}
we derive
$$
q_n\leq \frac{2\nu}{27}-\frac{\pi\nu n^2}{2|\Omega|}-\frac{\pi^2\nu_0 n^4}{|\Omega|^2}+\frac{C}{\nu}\sup_{t\in\mathbb{R}}\sup_{u_0\in\mathcal{A}(t)}\Big(\frac{1}{T}\int^t_{t-T}\|S(s,\tau)u_0\|d s\Big)
$$
and
$$
\hat{q}_n\leq \frac{2\nu}{27}-\frac{\pi^2\nu_0 n^4}{|\Omega|^2}-\frac{\pi\nu n^2}{2|\Omega|}+\frac{C}{\nu}q,
$$
where $q=\displaystyle{\limsup_{T\rightarrow +\infty}}\sup_{t\in\mathbb{R}}\sup_{u_0\in\mathcal{A}(t)}\frac{1}{T}\int^t_{t-T}\|S(s,\tau)u_0\|^2d s$.

From the estimate of equation, we have
\begin{eqnarray}
\frac{\nu}{2}\int^t_s\|u(r)\|^2dr+2\nu_0\int^t_s\|u(r)\|^4dr\leq |u_0(s)|^2+C\|f(r)\|_{L^2(s,t;V')}.\label{yan-15}
\end{eqnarray}
Setting $s=t-T$ in \eqref{yan-15}, using the bounded of solution in $H$, it follows
\begin{align*}
q & \leq \frac{C}{\nu}\lim_{T\rightarrow +\infty}\frac{|u_0(t-T)|^2}{T}+\frac{C}{\nu}\lim_{T\rightarrow +\infty}\frac{1}{T}\int^t_{t-T}\|f(r)\|^2_{V'}dr\nonumber\\
& \leq  \frac{C}{\nu}\lim_{T\rightarrow +\infty}\frac{1}{T}\int^t_{t-T}\|f(r)\|^2_{V'}dr.
\end{align*}
Defining $M=\displaystyle{\lim_{T\rightarrow +\infty}}\frac{1}{T}\int^t_{t-T}\|f(r)\|^2_{V'}dr
= \langle \|f(\cdot) \|^2_{V'} \rangle|_{\le t}$ (cf. \eqref{G}), then
$$
\hat{q}_n \leq \frac{2\nu}{27}-\frac{\pi\nu n^2}{2|\Omega|}-\frac{\pi^2\nu_0 n^4}{|\Omega|^2}+\frac{C}{\nu^2_0}M.
$$

{\bf Case 1:} If $\frac{\pi\nu n^2}{2|\Omega|}+\frac{\pi^2\nu_0 n^4}{|\Omega|^2}>\frac{2\nu}{27}+\frac{C}{\nu^2_0}M$, then by Lemma 4.19 in \cite{clr}, we have $\mbox{dim}_B(\mathcal{A}(t))\leq 3$ with $n=3$.

{\bf Case 2:} Denoting $\hat{M}=\displaystyle{\limsup_{T\rightarrow +\infty}}\frac{1}{T}\int^t_{t-T}|f(r)|^2dr$, otherwise, by the theory in \cite{clr},  we see the fractal dimension of pullback attractors satisfies
$$
\mbox{dim}(\mathcal{A}(t))\leq \frac{C|\Omega|^{\frac{1}{2}}}{\nu^2_0\lambda_1}\hat{M}+\frac{2\nu}{27}=\hat{C}G+\frac{2\nu}{27},
$$
where here $G=\frac{\hat{M}}{\nu^2_0\lambda_1}$. The proof is completed.   \qed

\begin{Remark} \label{rem-fractal}
(1)	From the proof above, we can see that for the particular case $\frac{9\pi\nu}{2|\Omega|}+\frac{81\pi^2\nu_0}{|\Omega|^2}>\frac{C}{\nu^2_0}M+\frac{2\nu}{27}$,
then the fractal dimension of pullback attractor is estimated by $dim_F(\mathcal{A}^H_{\mu}(t))\leq 3$, where
 $C$ is a constant which depends on the bounded domain and first eigenvalue of Laplace
operator.
	
(2) The fractal dimension of global attractor $\mathcal{A}$ of 2D autonomous Navier-Stokes equation \eqref{bmy1-a} in $H$
can be estimated as
$$
\mbox{dim}_F\mathcal{A}  \leq CG, \quad C \leq (\frac{2}{\pi})^{1/2}(\lambda_1|\Omega|)^{1/2},
$$
for Dirichlet boundary condition, and
$$
\mbox{dim}_F\mathcal{A}  \leq C G^{2/3}(1+\log G)^{1/3},
$$
for periodic boundary conditions, with $G=\frac{|f|^2}{\nu^2\lambda_1}$. See \cite{fmtt,te}.
These results are also true for pullback attractors $\mathcal{A}(\cdot)$ in some non-autonomous cases, cf. \cite{clr,llr}.
	
(3) For \eqref{bmy1-a} in 3D, Chepyzhov and Ilyin \cite{ci} gave the estimate of invariant sets $X_A$ in $V$ as $\mbox{dim}_HX_A\leq CG^3$ and $\mbox{dim}_F X_A\leq 2CG^3$. However, since the uniqueness of global weak solution for 3D equation is still an open problem, estimates on the fractal dimension of trajectory attractor is unknown.
More results concerned to fractal dimension of attractors, we can refer to \cite{ci,cf85,cf1,cft84,cft,cv3,ily1,l1,llr,lu,te}.
\end{Remark}

\subsection{Proof of Theorem \ref{th5-g}.} \label{se3.8}
The pullback attractors $\mathcal{A}$ in $H$ becomes a single trajectory for some special viscosity $\nu, \nu_0$.
Let $u(t), v(t)$ be two solutions of problem \eqref{bmy1-1} with initial data $u(\tau)=u_0$ and $v(\tau)=v_0$ respectively.
Denote $w(t)=u(t)-v(t)$ and assume $\|u(t)\| \geq \|v(t)\|$ (or else denote $w=v-u$), then we see that $w$ satisfies
\begin{equation}\label{bmy1-0}
\begin{cases}
w_t-\Big[\nu+\nu_0(\|u(t)\|^2-\|v(t)\|^2)\Big]\Delta w +(u\cdot \nabla)w+(w\cdot\nabla)v=0,& \\
\nabla\cdot w=0,&\\
w|_{\partial\Omega}=0,& \\
w(x,\tau)=u_0-v_0. &
\end{cases}
\end{equation}
Multiplying \eqref{bmy1-0} with $w$, using Poincar\'{e}'s inequality and the property of $b(\cdot,\cdot,\cdot)$, it follows
\begin{align*}
\frac{1}{2}\frac{d|w|^2}{d t}+[\nu+\nu_0(\|u\|^2-\|v\|^2)]\|w\|^2\leq& |w|^{\frac{1}{2}}\|v\||w|^{\frac{3}{2}}\\
\leq& \frac{c}{\nu}\|w\|^2\|v\|^4+\frac{\nu}{2}\|w\|^2,
\end{align*}
and hence,
$$
\frac{d|w|^2}{dt}\leq \Big[\frac{c}{\nu}\|v\|^4+2\nu_0\lambda_1\|v\|^2-\nu\lambda_1-2\nu_0\lambda_1\|u\|^2\Big]\|w\|^2.
$$
If $u_0$ and $v_0$ fixed, we let $\tau$ goes to $-\infty$, it follows $u(t)=v(t)$, which means the pullback attractors is a point provided that
$$
\frac{c}{\nu}\|v\|^4+2\nu_0\lambda_1\|v\|^2-\nu\lambda_1-2\nu_0\lambda_1\|u\|^2<0.
$$
A sufficient but may be not optimal condition is
\begin{equation}
\frac{c}{\nu}\|v\|^4+2\nu_0\lambda_1\|v\|^2<\nu\lambda_1.\label{an-1}
\end{equation}
From the procedure of pullback absorbing set, we see that
$$
\frac{d}{d t}|v|^2+2(\nu+\nu_0\|v\|^2)\|v\|^2\leq \frac{2|f|^2}{\nu\lambda_1}+\nu\|v\|^2
$$
and
$$
\nu\int^t_{s}\|v\|^2d r+2\nu_0\int^t_s\|v\|^4d r\leq (|v(t)|^2-|v(s)|^2)+\frac{2}{\nu\lambda_1}\int^t_s|f(r)|^2d r,
$$
which implies
\begin{eqnarray}
\langle\|v\|^2\rangle|_{\leq t}\leq \frac{\langle2|f|^2\rangle|_{\leq t}}{\nu^2\lambda_1},\ \ \langle\|v\|^4\rangle|_{\leq t}\leq \frac{\langle|f|^2\rangle|_{\leq t}}{\nu\nu_0\lambda_1}.\label{an-2}
\end{eqnarray}
Combining \eqref{an-1} and \eqref{an-2}, yields
$$
\frac{c\langle|f|^2\rangle|_{\leq t}}{\nu^2\nu_0\lambda_1}+\frac{4\nu_0\lambda_1\langle|f|^2\rangle|_{\leq t}}{\nu^2\lambda_1}<\nu\lambda_1.
$$
From \eqref{Gg} we see that
${\rm G^g}(t)^2=\frac{\langle|f|^2\rangle|_{\leq t}}{\nu^4_0\lambda^2_1}$.
Then we can derive a sufficient condition for pullback attractors which is nontrivial when
$$
{\rm G^g}(t) < \sqrt{\frac{\nu_0}{c\nu+4\nu_0^2\nu\lambda_1}}.
$$
This ends the proof. \qed

%%%%%%%%%%%%%%%%%%%%%%%%%%%%%%%%%%%%%%%%%%%%%%%%%%%%%%%%

\subsection{Proof of Theorem \ref{re}.} \label{pro}
The pullback attractor $\mathcal{A}$ is the same considered in \cite{cf94,cf97}. Noting that
\begin{equation}
\mathcal{A}_{CDF}(t)=\overline{\bigcup_{B\ bounded\ in\ H}\Lambda(B,t)}^H,
\end{equation}
where $\Lambda(B,t)=\displaystyle{\bigcap_{s\leq t}}\ \overline{\bigcup_{t\leq s}U(t,\tau)B(\tau)}^H$, and since the universes $\mathcal{D}_F,\ \mathcal{D}_{\mu}$ and $\mathcal{D}_{\mu_0}$ in $H$ is no need to be bounded, $\mathcal{D}_{\mu}$ is arbitrary, it follows that
\begin{equation}
B\subset \mathcal{D}_F\subset \mathcal{D}_{\mu}\subset \mathcal{D}_{\mu_0}.
\end{equation}
Using the structure of pullback attractors $\mathcal{A}^H_{\mu}$ in Theorem \ref{th5.1}, i.e., the property of pullback-$\omega$ limit set
\begin{eqnarray}
\mathcal{A}^H_{\mu}=\displaystyle{\bigcap_{T\leq t}}\ \overline{\displaystyle{\bigcup_{s\leq T}}S(t,s)D_{\mu}(s)}^H,
\end{eqnarray}
we conclude that $\mathcal{A}_{CDF}(t)$ is included in other pullback attractors and
\begin{eqnarray}
\mathcal{A}_F(t)\subset \mathcal{A}^H_{\mu}\subset \mathcal{A}^H_{\mu_0}.
\end{eqnarray}
The similar result also holds in $D(A^{\frac{1}{2}})$,
Which implies (a) and (b).

From the theory in \cite{gmr,mr}, if the union of universes or pullback absorbing sets in uniformly bounded, then (c)-(f) is true. The proof has been completed. \qed

\subsection{Upper semi-continuity theory of pullback attractors}\label{sub3.3}
Consider the non-autonomous system with perturbed external force
\begin{equation}
\frac{\partial u}{\partial t}=\hat{A}_{F}u+\varepsilon F(t,x),\label{yang-1}\end{equation} our goal of this section is to show the relationship between pullback attractors ${\mathcal
	A}_{\varepsilon}=\{A_{\varepsilon}(t)\}_{t\in {\mathbb R}}$ and global
attractor $\mathcal {A}$ for \eqref{yang-1} with the cases $\varepsilon>0$ and
$\varepsilon=0$ respectively. The upper semi-continuity of attractors was investigated firstly by Hale and Raugel \cite{jkh} in 1988, then many mathematicians extended the theory to pullback attractor and random attractors for processes (cocycle), see Caraballo, Langa and Robinson \cite{ca}, Carvalho, Langa and Robinson \cite{clr}, Kloeden and Stonier \cite{klo}, Wang and Qin \cite{wq2010} and references therein.

In what follows, we will show the upper semi-continuity of pullback attractors  with respect to the parameter $\varepsilon\in (0,\varepsilon_{0}]$ for the evolutionary process
$U_{\varepsilon}(\cdot,\cdot)$ of \eqref{yang-1}.

For each $\tau\leq t\in {\mathbb R}$ and $x\in X$, we
assume
\begin{equation} \label{H1}
\lim \limits_{\varepsilon \rightarrow 0} {\rm dist}_{X}(U_{\varepsilon}(t,t-\tau)x, S(t-\tau)x)=0
\end{equation}
holds uniformly on bounded sets of $X$.

\begin{Definition}
	(See \cite{clr}) Let $X$ be a Banach space, $\Lambda$ be a metric space and ${A_{\lambda}} (\lambda\in\Lambda)$ be a family of subsets of
$X$. We say that the family of pullback attractors $A_{\lambda}$ is upper semi-continuous
as $\lambda\rightarrow \lambda_0$ if
$$
\lim_{\lambda\rightarrow \lambda_0} {\rm dist}_X(A_{\lambda},A_{\lambda_0})=0.
$$
\end{Definition}

\begin{Theorem}\label{th3} (See \cite{ca}) Assume that \eqref{H1} holds and there exist pullback attractors
${\mathcal A}_{\varepsilon}=\{A_{\varepsilon}(t)\}_{t\in {\mathbb R}}$ for all
$\varepsilon \in (0,\varepsilon_{0}]$. If there exists a compact set $K\subset X$, such that
\begin{equation} \label{H2}
\lim\limits_{\varepsilon\rightarrow 0}
{\rm dist}_{X}(A_{\varepsilon}(t), K)=0, \quad t\in \mathbb{R}.
\end{equation}
Then ${\mathcal A}_{\varepsilon}$ are upper semi-continuous to $\mathcal {A}$, i.e.,
$$
\lim\limits_{\varepsilon\rightarrow 0} {\rm dist}_{X}(A_{\varepsilon}(t),{\mathcal {A}})=0, \quad t\in \mathbb{R}.
$$
\end{Theorem}

In the sequel we present a procedure to verify \eqref{H2}.

\begin{Theorem}\label{th1}
	(See \cite{wq2010}) Assume the family of sets ${\mathcal B}=\{B(t)\}_{t\in{\mathbb R}}$ is
	pullback absorbing for the process $U(\cdot,\cdot)$, ${\mathcal
		K}_{\varepsilon}=\{K_{\varepsilon}(t)\}_{t\in {\mathbb R}}$ is a
	family of compact sets in $X$ for each
	$\varepsilon\in (0,\varepsilon_{0}]$. Suppose the decomposition
	$U_{\varepsilon}(\cdot,\cdot)=U_{1,\varepsilon}(\cdot,\cdot)+U_{2,\varepsilon}(\cdot,\cdot):{\mathbb
		R}\times{\mathbb R}\times X\rightarrow X$ satisfies
	
	(i) for any $t\in {\mathbb R}$ and $\varepsilon\in
	(0,\varepsilon_{0}]$,
$$
\parallel U_{1,\varepsilon}(t,t-\tau)x_{t-\tau}\parallel_{X}
\leq \Phi(t,\tau),\quad \forall \, x_{t-\tau} \in B(t-\tau) , \quad \tau>0,
$$
where
	$\Phi(\cdot,\cdot):{\mathbb R}\times{\mathbb R}\rightarrow {\mathbb
		R}^{+}$ satisfies $\lim\limits_{\tau\rightarrow
		+\infty}\Phi(t,\tau)=0$ for each $t \in \mathbb{R}$.
	
	(ii)  for any $t\in{\mathbb R}$ and $T\geq 0$, $\displaystyle{\bigcup_{0\leq
			\tau\leq T}}U_{2,\varepsilon}(t,t-\tau)B(t-\tau)$ is bounded, and for
	any $t\in {\mathbb R}$, there exists a time $T_{{\mathcal B}}(t)>0$,
	which is independent of $\varepsilon$, such that
$$
U_{2,\varepsilon}(t,t-\tau)B(t-\tau)\subset
	K_{\varepsilon}(t),\quad \forall \, \tau \geq T_{{\mathcal B}}(t), \quad
	\varepsilon\in (0,\varepsilon_{0}],
$$
and there exists a compact set
	$K\subset X$, such that
$$
\lim\limits_{\varepsilon\rightarrow 0} {\rm dist}_{X}(K_{\varepsilon}(t),K)=0,
\quad t \in \mathbb{R} .
$$
Then (a) for each $\varepsilon\in (0,\varepsilon_{0}]$, the system \eqref{yang-1}
	possesses a family of pullback attractors ${\mathcal A}_{\varepsilon} = \{A_{\varepsilon}(t)\}_{t\in {\mathbb R}}$,
(b) condition \eqref{H2} holds and hence $\mathcal{A}_{\varepsilon}$ is upper semi-continuous at $0^{+}$.
\end{Theorem}

\begin{Remark}
	In order to obtain the upper semi-continuity of attractors of the system \eqref{yang-1}, the weak solution must have the same initial data, i.e., every trajectory should begin at the same point.
\end{Remark}

\subsection{Proof of Theorem \ref{th5.10}} \label{sub5.2}
Using the theory in Section \ref{sub3.3}, we shall use the decomposition of process to estimate the linear equation with non-homogeneous initial data and nonlinear equation with homogeneous initial data, i.e.,
the solution
$u_{\varepsilon}(t)=U_{\varepsilon}(t,\tau)u_{\tau}$ of perturbed problem \eqref{yang-2} with and
$f(x,t)=\varepsilon h(x,t)$ and
initial data $u_{\tau}\in H$ can be decomposed as
$$
u_{\varepsilon}=S_{\varepsilon}(t,\tau)u_{\tau}=S_{1,\varepsilon}(t,\tau)u_{\tau}
+ S_{2,\varepsilon}(t,\tau)u_{\tau},
$$
where $ S_{1,\varepsilon}(t,\tau)u_{\tau}=v(t)$ and
$S_{2,\varepsilon}(t,\tau)u_{\tau}=w(t)$ solve the
problems
\begin{equation} \left\{\aligned &v_{t}+\nu Av+\mathbb{A}v=0,\\
&v(x,t)|_{\partial\Omega}=0, \\
&v(\tau,x)=u_0(x), \endaligned\right.\label{f17}
\end{equation}
and
\begin{equation} \left\{\aligned &w_{t}+\nu A w+\mathbb{A}w=-B(u,u)+\varepsilon h(x,t),\\
&w(x,t)|_{\partial\Omega}=0, \\
&w(\tau,x)=0, \endaligned\right.\label{f18}
\end{equation}
respectively.

\begin{Lemma} \label{le5.11-1} Let
	 $R_{\eta}=\{r:\mathbb{R}\rightarrow(0,+\infty)|\displaystyle{\lim_{\xi\rightarrow-\infty}}e^{\eta
		\xi}r^2(\xi)=0\}$ and denote by ${\mathcal{D}}_{\eta}$ the class of
	families $\hat{D}=\{D(t):t\in \mathbb{R}\}\subset {\mathcal{D}}(H)$ as universe such that
	$D(t)\subset \bar{B}(0,r_{\hat{D}}(t))$,
	where $\bar{B}(0,r_{\hat{D}}(t))$ is the closed ball in $H$
	centered at zero with radius $r_{\hat{D}}(t)$.
	Suppose that $u_0\in H$, the external force $h(x,t)\in L^2(\mathbb{R};H)$ satisfies \eqref{ya-1}. Then for any bounded set $B\subset H$
	and any fixed $t\in {\mathbb R}$, there exists a time $T(B,t)>0$, such that
$$
	\parallel S_{\varepsilon}(t,t-\tau)u_{t-\tau}\parallel^{2}_H\leq
	R^2_{\varepsilon}(t) \quad \forall \, \tau\geq T(B,t), \; u_{t-\tau}\in B,
$$
where
	$R^2_{\varepsilon}(t)=\frac{2C\varepsilon}{\nu} e^{-\eta t} \int_{-\infty}^{t}e^{\eta s}|h(s)|^{2}d s$.
	
	Moreover, setting $B_{\varepsilon}(t)=\{u_{\varepsilon} \in H \, | \,
	|u_{\varepsilon}|^2\leq R^2_{\varepsilon}(t)\}$,
	then ${\mathcal B}_{\varepsilon}=\{B_{\varepsilon}(t)\}_{t\in {\mathbb{R}}}\in \mathcal{D}_{\eta}$ is
	the family of pullback absorbing sets in $H$, i.e.,
\begin{equation}
	\lim\limits_{t\rightarrow -\infty}e^{\eta t}R_{\varepsilon}(t)=0
	\quad \forall ~\varepsilon>0. \label{z4}
\end{equation}
\end{Lemma}

\noindent{\bf Proof.}
	Let $t\in \mathbb{R}$ be fixed, then for any $\tau\in \mathbb{R}$ and $u_0\in H$, we denote
$$
u_{\varepsilon}(r)=u(r;t-\tau,u_0)=u_{\varepsilon}(r-t+\tau,t-\tau,u_0) \quad \forall \, r \geq t-\tau.
$$
	Multiplying perturbed problem \eqref{yang-2} ($f(t)=\varepsilon h(x,t)$) with $e^{\eta t}u_{\varepsilon}$ ($\eta$ will be determined later), noting that $(B(u_{\varepsilon},u_{\varepsilon}),u_{\varepsilon})=0$, we derive that
	\begin{align}
	\frac{d}{d t} \Big(e^{\eta t}|u_{\varepsilon}(t)|^2\Big) & +2\nu e^{\eta	 t}\|u_{\varepsilon}(t)\|^2+2\nu_0 e^{\eta	 t}\|u\|^2\|u_{\varepsilon}(t)\|^2\nonumber\\
	& =  \eta e^{\eta t}|u_{\varepsilon}(t)|^2+2e^{\eta t}(\varepsilon h(t), u_{\varepsilon}(t))\nonumber\\
	& \leq  \eta e^{\eta t}|u_{\varepsilon}(t)|^2+\nu e^{\eta t}\|u_{\varepsilon}(t)\|^2+\frac{C\varepsilon}{\nu}e^{\eta t}|h(t)|^2,\label{f21}
	\end{align}
	holds for all $u_{\varepsilon}\in H$, then using the Poincar\'{e} inequality, choosing $\eta=\frac{\nu\lambda_1}{2}$ and neglecting the third term in \eqref{f21}, we have
$$
\frac{d}{d t}\Big(e^{\eta t}|u_{\varepsilon}(t)|^2\Big)+\frac{\nu\lambda_1}{2} e^{\eta
		t}|u_{\varepsilon}(t)|^2\leq \frac{C\varepsilon}{\nu}e^{\eta
		t}|h(t)|^2,
$$
which implies
$$ |u_{\varepsilon}(t)|^2\leq e^{-\eta(t-\tau)}\| u_{0}\|^2+\frac{C\varepsilon}{\nu}\int^t_{\tau}
 e^{-\eta(t-\xi)}|h(\xi)|^2d \xi
$$
for all $\tau\in \mathbb{R}$.
	
Let $\hat{D}\in {\mathcal{D}}_{\eta}$ be given above, then for any $u_{0}\in D(\tau)$ and $t\geq\tau$,
it yields
$$
|S_{\varepsilon}(t,t-\tau)u_{t-\tau}|^2\leq e^{-\eta(t-\tau)}r^2_{\hat{D}}+\frac{C\varepsilon}{\nu}
	 \int^t_{-\infty}e^{-\eta(t-\xi)}|h(\xi)|^2d\xi.
$$
Setting $e^{-\eta(t-\tau)}r^2_{\hat{D}}\leq\frac{C\varepsilon}{\nu}\int^t_{-\infty}e^{-\eta(t-\xi)}|h(\xi)|^2d\xi$,
then for fixed $t\in \mathbb{R}$,
we denote $R_{\varepsilon}(t)>0$ as
$$
(R_{\varepsilon}(t))^2=\frac{2C\varepsilon}{\nu}
	 \int^t_{-\infty}e^{-\eta(t-\xi)}|h(\xi)|^2d\xi.\label{f25}
$$
Considering the family of closed balls $\hat{B}_{\varepsilon}$ for any fixed $t\geq \tau$ in $H$ defined by
$$
B_{\varepsilon}(t)=\{u_{\varepsilon}\in H \, | \, |u_{\varepsilon}|^2\leq
	2R^2_{\varepsilon}(t)\},
$$
it is easy to check that $\mathcal{B}_{\varepsilon}(t)\in {\mathcal{D}}_{\eta}$ and hence
	$\mathcal{B}_{\eta}(t)$ is the family of ${\mathcal{D}}_{\eta}$-pullback absorbing
	sets for the process $\{S_{\varepsilon}(t,t-\tau)\}$. \qed

\begin{Lemma} \label{le5.11} Let $R_{\varepsilon}(t)$, $B_{\varepsilon}(t)$ are defined
in Lemma \ref{le5.11-1}, then for any $t\geq \tau\in {\mathbb R}$, the solution $v(t)=S_{1,\varepsilon}(t,t-\tau)u(t-\tau)$ of 	 \eqref{f17} satisfies
\begin{align}
|S_{1,\varepsilon}(t,t-\tau)u_{t-\tau}|^{2} & \leq e^{-2\nu\lambda_1\tau}R^2_{\varepsilon}(t-\tau),\nonumber \\
\int^t_{t-\tau}\|v(s)\|^2d s & \leq J_{\varepsilon}(t)\label{a25-a}
\end{align}
	for all $\tau\in \mathbb{R}$ and $u_{t-\tau}\in B_{\varepsilon}(t-\tau)$, where $J_{\varepsilon}(t)$ is dependent on $\tau, \ R^2_{\varepsilon}(t-\tau),\ \nu$ and $\lambda_1$.
\end{Lemma}

\noindent{\bf Proof.} Multiplying \eqref{f17} with $v$ and integrating by part over $\Omega$, we obtain
	\begin{align} \frac{1}{2}\frac{d}{d t}|
	v(t)|^2+\nu\|v(t)\|^{2}+\nu_0\|u\|^2\|v\|^2\leq 0,\label{ya-3}\end{align}
	here we use $(B(v,v),v)=0$. By using the Poincar\'{e} inequality, neglecting the third term in \eqref{ya-3}, it yields
	\begin{equation} \frac{d}{d t}|v(t)|^2+2\nu\lambda_1|v(t)|^{2}\leq 0.\label{a25-1}\end{equation}
	Applying Gronwall's inequality to \eqref{a25-1} from $t-\tau$ to $t$, we get
	\begin{equation}|S_{1,\varepsilon}(t,t-\tau)v_{t-\tau}|^{2}\leq| v_{t-\tau}|^2e^{-2\nu\lambda_1\tau}\leq
	e^{-2\nu\lambda_1\tau}R^2_{\varepsilon}(t-\tau), \forall\ t\geq\tau . \label{a261}
	\end{equation}
	Estimate \eqref{a25-a} is the direct result of \eqref{a261}. This completes the proof. \qed

\begin{Lemma}\label{le5.13} Let the family of pullback absorbing sets ${\mathcal B}_{\varepsilon}(t)=\{B_{\varepsilon}(t)\}_{t\in {\mathbb R}}$ be given by Lemma \ref{le5.11-1} and \eqref{z4} holds. Then,  for any fixed $t\geq \tau\in {\mathbb R}$, there exist a time
	$T_{\varepsilon}(t,{\mathcal B})>0$ and a function $I_{\varepsilon}(t)>0$, such that
	the solution $S_{2,\varepsilon}(t,\tau)u_{\tau}=w(t)$ of \eqref{f18} satisfies
$$
 \|S_{2,\varepsilon}(t,t-\tau)u_{t-\tau}\|^{2}_{D(A^{\frac{1}{2}})}\leq
	I_{\varepsilon}(t),
$$
for all  $\tau\geq
	T_{\varepsilon}(t,{\mathcal B}) $ and any $u_{t-\tau}\in
	B_{\varepsilon}(t-\tau)$.
\end{Lemma}

\noindent{\bf Proof.} Taking inner product of \eqref{f18} with $A w(t)$ in $H$, integrating by parts over $\Omega$, we derive
\begin{equation} \label{a28}
\frac{1}{2}\frac{d}{d t}|A^{\frac{1}{2}}w(t)|^2+\nu|A w(t)|^2+\nu_0\|u\|^2|A w(t)|^2
	 =-b(u,u,A w)+\varepsilon\langle h(t), Aw\rangle.
\end{equation}	
	By the property of trilinear operator $b(\cdot,\cdot,\cdot)$ and Young's inequality, we obtain
\begin{align*}
 |b(u,u,A w)|
 & \leq |u|^{1/4}\|u\|^{3/4}\|u\|^{1/4}|Au|^{3/4}|A w| \\
 & \leq C|u|^2+C|Au|^2+\nu_0\|u\|^2|A w|^2 ,
\end{align*}
	and
\begin{equation}
\langle\varepsilon h(t),A w\rangle\leq \frac{\nu}{2}|A w(t)|^2+\frac{C\varepsilon^2}{\nu}|h(t)|^{2}.\label{a30}
\end{equation}
Hence combining \eqref{a28}-\eqref{a30}, we derive
\begin{equation}\frac{d}{d t}|A^{\frac{1}{2}}w(t)|^{2}+\nu|Aw(t)|^2\leq C|u|^2+C|Au|^2+\frac{C\varepsilon^2}{\nu}|h(t)|^{2}.\label{a31}
\end{equation}
Applying the Gronwall inequality to \eqref{a31} from $t-\tau$ to $t$, using Lemma \ref{le5.11}, we conclude that
$$
|A^{\frac{1}{2}}w(t)|^{2}\leq I_{\varepsilon}(t)=I_{\varepsilon}(t,\tau,R_{\varepsilon}(t-\tau),
J_{\varepsilon},\nu,\nu_0,\int^t_{-\infty}e^{\beta s}\|h(s)\|^2d s)
$$
	for all $t\geq\tau$. This achieve the proof of desired lemma. \qed

\begin{Lemma} \label{le5.14}For any fixed $t\geq\tau\in {\mathbb R}$, if $u_{0}$ takes its value in some bounded set, then the solution
	$u_{\varepsilon}(t)=S_{\varepsilon}(t,t-\tau)u_{0}$ of perturbed non-autonomous problem \eqref{bmy1-1},
with $f(x,t)=\varepsilon h(x,t)$, converges to the solution
	$u(t)=S(t)u_{0}$ of the autonomous problem with
	$f=0$ uniformly in $H$ as $\varepsilon \to 0^{+}$. This means
\begin{eqnarray} \lim\limits_{\varepsilon\rightarrow
		0^{+}}\sup\limits_{u_{0}\in B}\|u_{\varepsilon}(t)-u(t)\|_{H}=0,\label{a33}
\end{eqnarray}
	where $B$ is a bounded subset in $H$.
\end{Lemma}

\noindent{\bf Proof.}
Denoting
$$
y^{\varepsilon}(t)=u_{\varepsilon}(t)-u(t),
$$
we can verify that $y^{\varepsilon}(t)$ satisfies the problem
\begin{equation} \left\{\aligned &\frac{d y^{\epsilon}}{d t}+\nu Ay^{\varepsilon}+\mathbb{A}y^{\varepsilon}=-B(u_{\varepsilon},u_{\varepsilon})+B(u,u)+\varepsilon h(t,x),\\
	&y_{\varepsilon}|_{\partial\Omega}=0, \\
	&y_{\varepsilon}|_{t=\tau}=(u_{\varepsilon})_{\tau}-u_{\tau}=0.\\
	\endaligned\right.\label{a35}
\end{equation}
Multiplying \eqref{a35} by $y^{\varepsilon}(t)$, using the property of $b(\cdot,\cdot,\cdot)$, we have
\begin{align}
\frac{1}{2}\frac{d}{d t}|y^{\varepsilon}|^{2} & + \nu \|y^{\varepsilon}\|^2+\nu_0\|u\|^2\| y^{\varepsilon}\|^{2} \nonumber \\
& = \langle B(u,u)-B(u_{\varepsilon},u_{\varepsilon}),y^{\varepsilon} \rangle
+ \langle\varepsilon h(t),y^{\varepsilon}\rangle \nonumber \\
& \leq |\langle B(u,u)-B(u_{\varepsilon},u_{\varepsilon}),y^{\varepsilon})\rangle|
+\frac{\nu}{2}\|y^{\varepsilon}(t)\|^2+\frac{C\varepsilon^2}{\nu}|h(t)|^{2}.\label{ya-0}
\end{align}
	By Young's inequality, noting that $b(u_{\varepsilon},y^{\varepsilon}, y^{\varepsilon})=0$, we get
\begin{align*}
|\langle B(u,u)-B(u_{\varepsilon},u_{\varepsilon}),y^{\varepsilon})\rangle|
& = |b(y^{\varepsilon},u,y^{\varepsilon})| \\
& \leq \frac{C}{\nu}\|u\|^4|y^{\varepsilon}|^2+\nu\|y^{\varepsilon}\|^{2}.
\end{align*}
	Hence, neglecting the third term in \eqref{ya-0}, it follows
\begin{equation}\frac{d}{d t}|y^{\varepsilon}|^{2}\leq \frac{C}{\nu}|y^{\varepsilon}|^{2}\|u\|^{4}+\frac{C\varepsilon^2}{\nu}|h(t)|^{2}.\label{ya-6}
\end{equation}
Using Lemmas \ref{le5.11-1} to \ref{le5.14} and \eqref{ya-1}, noting that $h\in L^{2}_{loc}({\mathbb R},H)$,
using the Gronwall inequality to \eqref{ya-6}, we conclude
\begin{align*}
|y^{\varepsilon}|^{2}
& \leq \frac{C\varepsilon^2}{\nu}e^{\frac{C}{\nu}\|u\|^4_{L^4(\tau,T;V)}}\int_{t-\tau}^{t} |h(s)|^{2} \, ds \\
& \leq \varepsilon C_{\tau,t} \to 0 ,
\end{align*}
as $\varepsilon\rightarrow 0^+$, which implies \eqref{a33}. This ends the proof. \qed

\medskip

\noindent{\bf Proof of Theorem 3.15.} 	
Now the proof of the upper semi-continuity of pullback attractors in $H$ follows from Lemma \ref{le5.14}. \qed	

\section{Conclusion and further research}
From the discussion in this paper, we can see that the 3D Navier-Stokes equation with nonlinear viscosity \eqref{bmy1-1} has better dissipative property than the
classical 3D model \eqref{bmy1-a}. The disadvantage is that \eqref{bmy1-1} does not satisfy the Stokes principle. On the other hand, since the well-posedness of 3D Navier-Stokes equation is still an open problem, one could study the long-time dynamics of a class of physically justified Ladyzhenskaya models \eqref{bmy1-b} that satisfy the Stokes principle and are also well-posed. Moreover, the upper semi-continuity of pullback attractors to trajectory attractors of \eqref{bmy1-1}
as $\nu_0$ goes to $0$ is still an unsolved problem.

\medskip

\paragraph{Acknowledgements} This work was initiated when X. Yang was a long term visitor at ICMC-USP, Brazil, from May 2015 to June 2016, supported by FAPESP (Grant No. 2014/17080-0). He was also partially supported by NSFC of China (Grant No. 11726626). Y. Lu was partially supported by NSF (Grant No. 1601127). They
were also supported by the Key Project of Science and Technology of Henan Province (Grant No. 182102410069, 172102210342, 17A120003). B. Feng was supported by NSFC (Grant No. 11701465).
T. F. Ma was partially supported by CNPq (Grant No. 310041/2015-5).
The authors thank Professors Chunyou Sun (Lanzhou University) and Yonghai Wang (Donghua University)
for fruitful discussion on this subject.


\end{document}